\documentclass{amsart}

\usepackage[utf8]{inputenc}
\usepackage{amsmath,amsfonts,amsthm,amssymb,cite,xcolor} %blkarray?
\usepackage{hyperref} % For preprint URLs
\usepackage[capitalise, noabbrev, nameinlink]{cleveref}
\hypersetup{pdfborder={0 0 0}}

\usepackage{graphicx}
\usepackage{enumerate}

\usepackage[english]{babel}

\usepackage[english=usenglishmax]{hyphsubst}

\usepackage{tikz}
\usetikzlibrary{calc}

\copyrightinfo{2026}{H. Tracy Hall}

\usepackage{thmtools}
\usepackage{thm-restate}

\newtheorem{theorem}{Theorem}[section]
\newtheorem{lemma}[theorem]{Lemma}

\newtheorem{corollary}[theorem]{Corollary}
\newtheorem{observation}[theorem]{Observation}
\newtheorem{proposition}[theorem]{Proposition}

\theoremstyle{definition}
\newtheorem{definition}[theorem]{Definition}
\newtheorem{example}[theorem]{Example}

\theoremstyle{remark}
\newtheorem{remark}[theorem]{Remark}

\theoremstyle{plain} % For remaining one-time conjectures

\numberwithin{equation}{section}

\newcommand{\psd}{_{+}}

\newcommand{\dagg}[1]{#1^\top\mspace{-2mu}}

\newcommand{\ad}{d}  % ambient dimension

\DeclareMathOperator{\M}{M}

\DeclareMathOperator{\Mpu}{M\psd^u}

\newcommand{\gdel}{\ensuremath{{}^g\mspace{-4mu}\Delta}} % 3mu gives better kerning?

\newcommand{\st}{: \ }

\begin{document}

\title[Delta Theorem]{The Delta Theorem: A dimension bound for faithful orthogonal graph representations}

\author{H. Tracy Hall}
\email{h.tracy@gmail.com}
\thanks{Parts of this research were conducted for the American Institute of Mathematics
while visiting Brigham Young University.
The author thanks AIM, BYU, and the National Science Foundation for their support.}

\subjclass[2010]{Primary }

\date{}

\dedicatory{}

\begin{abstract}
In 1987 Hiroshi Maehara conjectured that every graph can be represented
by a collection of vectors that are considered adjacent when not orthogonal
(a
\emph{faithful orthogonal representation} of the graph)
in dimension complementary to the minimum degree of the graph.
Without settling the conjecture,
L\`asl\`o~Lov\`asz, Michael~Saks, and Alexander~Schrijver showed
that dimension complementary to vertex connectivity both suffices and is best
possible under the additional assumption of
general position.
The present work proves the conjecture of Maehara as well as related conjectures
that have arisen independently in combinatorial matrix theory,
the strongest of which, $\delta(G) \le \nu(G)$,
is that minimum vertex degree gives a lower bound for the maximum nullity of a positive definite matrix with pattern $G$ that has the Strong Arnold Property.
Nullity questions, with various matrix restrictions, are an important subcase of the Inverse Eigenvalue Problem for a Graph (IEPG).

The name \emph{greedegree} is introduced for the largest possible final degree of a
maximum cardinality search (MSC) ordering, which is to say an ordering that
greedily maximizes adjacencies to previous chosen vertices.
The name \emph{upper-zero generic} is introduced to describe
symmetric matrices with nonzero diagonal
such that the zeros above the diagonal in any column belong to an independent
set of rows.
It is shown that some faithful orthogonal graph representation whose
Gram matrix is upper-zero generic always exists in dimension complementary
to the greedegree of the graph.
\end{abstract}



% \linenumbers

\maketitle

\section{Introduction}
\label{sec:intro}
Orthogonal representations were introduced by L\'aszl\'o L\'ovasz in 1979 \cite{L1979}
for studying the Shannon capacity of a graph, and have since also found
application for example in quantum information theory \cite{AL2001}.
For a general graph $G$ on $n$ vertices it is a difficult question, with a wide literature of partial results,
to determine the minimum possible dimension
for a faithful orthogonal representation of $G$.
The present work defines a class of vertex orderings
and associated target dimensions for which it is shown that
an existing probabilistic construction will succeed with probability $1$.
The resulting dimension bound verifies long-standing conjectures
of Hiroshi Maehara and others.

Fix a graph $G$ on $n$ vertices.
A representation of the vertices of $G$ by vectors in Euclidean space
is called an \emph{orthogonal representation of $G$}
if non-adjacency of a vertex pair implies orthogonality of the representing
vectors,
and an orthogonal representation is called \emph{faithful}
when the reverse implication also holds.
The matrix of pairwise dot products is the \emph{Gram matrix} of the representation,
and representations with equal Gram matrices are considered equivalent,
even though the dimension of the ambient space into which
the vertices are represented may differ.
The dimension of the span of the representing vectors is the
same for equivalent representations, and is equal to the
rank of the Gram matrix.
Complementing the span dimension with respect to $n$ gives
the \emph{codimension} of an orthogonal representation,
equal to
the nullity of its Gram matrix.
The largest possible codimension of a faithful orthogonal representation
of $G$ is equal to the graph parameter $\mathrm{M}_+(G)$
from combinatorial matrix theory,
which together with $\mathrm{M}(G)$ constitutes an important sub-problem of the Inverse Eigenvalue Problem for a Graph (IEPG).
An overview of IEPG and related matrix-theoretic graph parameters is given in \cite[Chapter 46]{HLA2016}.

Codimension bounds in terms of the vertex connectivity $\kappa(G)$ are given
in work
\cite{LSS1989, LSS2000} by
L\`aszl\`o~Lov\`asz,
Michael Saks, and
Alexander Schrijver 
(hereafter LSS), who
introduce a probabilistic construction of orthogonal representations.
When the LSS construction is carried out in Euclidean space
of dimension $n - \kappa(G)$,
they show that it yields, with probability $1$,
a representation that is both faithful and in general position.
Conversely, they show that $\kappa(G)$ is the maximum possible codimension
for a representation
that is both faithful and in general position.
When general position is not required, less is known,
but they mention a conjecture involving the minimum
vertex degree $\delta(G)$.
The conjecture, attributed to Maehara in 1987,
is that every graph $G$ has a faithful orthogonal representation
of codimension $\delta(G)$, equivalent to the statement
$\delta(G) \le \mathrm{M}_+(G)$.
A related conjecture in combinatorial matrix theory known as the Delta Conjecture later arose
independently
\cite{AIM2006}. It is weaker in one regard:
Where Maehara's Conjecture expects a positive semidefinite matrix,
the Delta Conjecture expects only a real symmetric matrix.

The present results establish the conjectured bound.
The proven bound is stronger than Maehara's Conjecture in two ways, the
first of which is that
the codimension obtained is in general larger, as given by a graph parameter
that will be introduced and given the name ``maximal greedy final degree'' or
``greedegree'' of $G$, denoted $\gdel(G)$.
The second way in which the proven bound is stronger relates to
a condition variously called
the Strong Arnold Hypothesis or the Strong Arnold Property
(hereafter SAP)
that was introduced to combinatorial matrix theory by Yves Colin de Verdi\'ere
when defining the
graph parameter $\mu(G)$ \cite{CdV1990, CdV1991}.
The SAP can be thought of as a requirement for a matrix to be generic, but in a certain technical sense that is
subject to the constraints of a specified low rank and a specified sparse pattern, neither of which is a generic condition.
The strongest sense of ``generic'' that one might require of such a matrix is general position,
but by the converse results of LSS, general position is impossible
in codimension $\delta(G)$, and therefore also in codimension $\gdel(G)$,
in the case that the inequality
$\kappa(G) \le \delta(G)$ is strict.
Between the weakly generic requirement of SAP and the strongly generic requirement of general position,
a generic requirement of intermediate strength will be introduced, and matrices satisfying the
intermediate requirement will be called ``upper-zero generic''.
The present work will demonstrate that representations can always be found
in codimension $\gdel(G)$ whose Gram matrices are upper-zero generic, and are therefore also
generic in the sense required by SAP.
This suffices to prove not only Maehara's Conjecture but also
the Strong Delta Conjecture \cite{BBFHHSvv2012},
which expects a positive semidefinite matrix that satisfies SAP.
Of the three standing conjectures resolved by the present results, the weakest one, the Delta Conjecture,
is perhaps the best known in the combinatorial matrix theory literature, and has
led to some partial results \cite{BFHHSv2009, BFHRS2008}.

The method of proof starts by taking the
probabilistic construction of LSS
and completely parametrizing it, in a symmetry-preserving way that introduces redundant degrees of freedom.
Although this yields polynomials that are of a difficult size to work with when
written out in full,
the
polynomials are tamed by introducing a construction, called a
hanging garden diagram, that is
somewhat reminiscent of, but overall simpler than, the planar algebras of Vaughan Jones \cite{J1999}.
Faithfulness of the representation requires certain of these polynomials to take non-zero values
for generic parameter choices.
Using these diagrams, it becomes possible to verify,
for the required polynomials according to a certain
monomial term order,
that there is no cancellation in the leading term.
The method of careful
accounting to show non-cancellation in the leading term
of a multivariable polynomial
is a technique that has previously been used for example in
Daan Krammer's demonstration \cite{K2002} of the faithfulness
of the Lawrence--Krammer representation of the braid group
by matrices of polynomials.

The work is organized as follows:

\Cref{sec:back} provides basic definitions and necessary background.

\Cref{sec:statement} presents the LSS construction and its known results, and introduces a slight modification of LSS called uniform LSS, which may succeed either weakly or strongly.
It then introduces terminology, with examples, sufficient to state the main results.

\Cref{sec:hanging} presents the hanging garden construction and details many of its properties.

\Cref{sec:LSS-garden} considers those hanging gardens specifically that model uniform LSS, including their structure when $G$ has an MCS or greedy ordering,
leading up to the Unique Monomial Lemma that is then used to prove the main results.

\Cref{sec:more} concludes with some possible generalizations.

\Cref{sec:NP} establishes, by reduction from 3-SAT, that the new combinatorial bound of greedegree, $\gdel(G)$, is in general NP-hard to determine.

\subsection{Acknowledgments}
The author acknowledges the patience of the combinatorial matrix theory community during the long interval between the initial claim of these results and their eventual written form.

Parts of this research were conducted for the American Institute of Mathematics
while visiting Brigham Young University.
The author thanks AIM, BYU, and the National Science Foundation for their support.

Thanks are due to Raphael Loewy for an initial introduction to the problem, to Leslie Hogben and others for encouragement and opportunities for presentation and discussion, and to Virginia and Barry Wood for hosting the final push of writing.

\section{Background}
\label{sec:back}
\subsection{Notation and conventions}
As given by a fixed choice in the Introduction, $G$ is a simple undirected graph on $n$ vertices.
A fixed ordering of the vertices will also be assumed by labeling the vertex set $V$ as a list
$(v_1, \dots, v_n)$,
and certain constructions will depend on that ordering, as will an implied correspondence between
vertices and row and column indices in matrices and vectors.
It will sometimes be relevant to consider different vertex orderings for the same graph $G$,
which will be accomplished by the assumption that the vertices have been relabeled
so that $v_1$ is now the first vertex in the new ordering, $v_2$ the second, and so forth.

The edge $\{v_i, v_j\}$ is abbreviated $v_i v_j$.
Adjacency in $G$ is
indicated by $v_i \sim v_j$ and non-adjacency by $v_i \not \sim v_j$.
The graph with edge set complementary to $G$ is denoted by
$\overline{G}$.
The degree $\mathrm{deg}(v_i)$ of a vertex $v_i$ is the number of vertices $v_j$ with $v_i \sim v_j$,
and the minimum degree and maximum degree among vertices of $G$ are denoted $\delta(G)$ and $\Delta(G)$ respectively.
The vertex connectivity of a graph is denoted $\kappa(G)$.
The subgraph of $G$ induced by a subset $S \subseteq V(G)$
is denoted $G[S]$.

The identity matrix is denoted $I$.
%, or $I_k$ to specify an identity matrix of size $k \times k$.
Given a square matrix $A$, the principal submatrix of $A$ on the set of indices $S$ is denoted $A[S]$.
The $n$ vertices of $G$ are assumed to be represented as column vectors
in some ambient dimension $\ad$
by a function $\varphi: V \rightarrow \mathbb{R}^\ad$,
whose outputs are collected into a result matrix of size $\ad \times n$,
$R = [r_{ij}] = [\varphi(v_0)\dots\varphi(v_{n-1})].$
The transpose of a matrix $A$ is denoted $\dagg{A}$.
The Gram matrix of $\varphi$ is $\dagg{R}R$, which
will be named $T = [t_{ij}]$,
where $t_{ij}$ is the inner product $\dagg{\varphi(v_i)} \varphi(v_j)$.

\subsection{Orthogonal
and faithful orthogonal representations}
The map $\varphi$
is called an \emph{orthogonal representation} of $G$
if, for any distinct vertices $v_i$ and $v_j$ of $G$
with $v_i \not \sim v_j$,
the inner product $t_{ij}$ of their images under $\varphi$ is $0$.
(For any graph $G$, the zero map on $n$ vertices
provides a trivial example of an orthogonal representation of $G$.)
If in addition the length of every vector $\varphi(v_i)$ is $1$,
$\varphi$ is called an \emph{orthonormal representation,}
in which case $T$ is a correlation matrix.
Two representations of the same graph are considered equivalent
if and only if they have the same Gram matrix $T$,
or in other words if and only if their images are isometric.
The \emph{codimension} of the orthogonal representation is the nullity
of $T$, which is at least $n - \ad$, with equality in the case that the image of $\varphi$ spans $\mathbb{R}^\ad$.

An orthogonal or orthonormal representation is called \emph{faithful} when the implication
carries in both directions:
In a faithful orthogonal representation, the distinct pairs
$i$ and $j$ for which $t_{ij} = 0$ holds are exactly those for which $v_i \not \sim v_j$,
so that the vectors of the representation
completely determine $G$ by the pattern of nonzero
inner products in the off-diagonal entries of the Gram matrix $T$.

The goal for a given graph $G$ is to produce a faithful orthogonal representation
whose codimension is as large as possible.
This turns out to relate directly to the maximum semidefinite nullity problem
of combinatorial matrix theory.
\begin{definition}
\label{def:Mplus}
Let $G$ be a graph on $n$ vertices. The  \emph{maximum semidefinite nullity} of $G$,
denoted $\mathrm{M}_+(G)$, is the largest nullity among all matrices $A = [a_{ij}]$
such that
\begin{enumerate}
\item the matrix $A$ of size $n \times n$ is symmetric and positive semidefinite,
\item for %$0 \le i < j \le n - 1$,
$i \ne j$, $a_{ij} = 0$ holds whenever $v_i \not\sim v_j$ in $G$, and
\item for %$0 \le i < j \le n - 1$,
$i \ne j$, $a_{ij} \ne 0$ holds whenever $v_i \sim v_j$ in $G$.
\end{enumerate}
\end{definition}

\begin{observation}
The conditions on $A$ in \cref{def:Mplus} are equivalent to saying:
\begin{enumerate}
\item There is a Euclidean arrangement of $n$ vectors whose Gram matrix is $A$.
\item The vector arrangement is an orthogonal representation of $G$.
\item The orthogonal representation is faithful.
\end{enumerate}
Since codimension is given by the nullity of the Gram matrix,
it follows that the maximum possible codimension of a
faithful orthogonal representation of $G$ is precisely $\mathrm{M}_+(G)$.
\end{observation}

Maximum semidefinite nullity is a special case of maximum nullity.

\begin{definition}
\label{def:M}
Let $G$ be a graph on $n$ vertices. The  \emph{maximum nullity} of $G$,
denoted $\mathrm{M}(G)$, is the largest nullity among all matrices $A = [a_{ij}]$
such that
\begin{enumerate}
\item the matrix $A$ of size $n \times n$ is real symmetric,
\item for % $0 \le i < j \le n - 1$,
$i \ne j$, $a_{ij} = 0$ holds whenever $v_i \not\sim v_j$ in $G$, and
\item for %$0 \le i < j \le n - 1$,
$i \ne j$, $a_{ij} \ne 0$ holds whenever $v_i \sim v_j$ in $G$.
\end{enumerate}
\end{definition}

Equivalent to the maximum nullity problem, and perhaps more widely referenced,
is the minimum rank problem.
The \emph{minimum rank} $\mathrm{mr}(G)$ of a graph $G$ is defined in such a way that
it satisfies $\mathrm{mr}(G) + \mathrm{M}(G) = n$.

\begin{example}[The Laplacian representation]
With the goal eventually of constructing a faithful orthogonal representation of high codimension,
it is instructive to consider an easy construction that always
yields a faithful orthogonal representation of codimension at least one.
The \emph{Laplacian representation} $\varphi_\mathcal{L}$ is a map from $V$
to $\mathbb{R}^{|E|}$, where the edges of $G$ have been numbered
$1$ to $|E|$ and each edge has been given an orientation.
(Different orientations will lead to isometric representations, so the choice is arbitrary.)
The column vector $\varphi_\mathcal{L}(v_i)$ has entry $-1$ for each edge
adjacent to and oriented away from $v_i$, entry $1$ for each edge
adjacent to and oriented towards $v_i$, and entry $0$ otherwise.
Equivalently, in the result matrix $R$, a row $k$ corresponding
to an edge which has been oriented $(v_i, v_j)$ has only
two non-zero entries, $r_{ki} = -1$ and $r_{kj} = 1$.
If $v_i \not \sim v_j$, then the supports of $\varphi_\mathcal{L}(v_i)$ and
$\varphi_\mathcal{L}(v_j)$ are disjoint and the vectors are orthogonal;
if $v_i \sim v_j$ then the inner product of the vectors is $-1$.
The ambient dimension $d = |E|$ is typically much larger than $n$,
and so $n - \ad$ is not a useful lower bound for the codimension
in this case.
Notice however that the sum of the $n$ representing vectors is zero,
which means that $T$ has a non-trivial null vector and the codimension
is at least one. (In fact,
the codimension of the Laplacian representation is exactly the
number of connected components of $G$.)
\end{example}

\section{Statements of results}
\label{sec:statement}
This section will review a series of conjectures related to codimensions of
orthogonal representations,
discuss the LSS probabilistic construction of orthogonal representations and a generalization called uniform LSS,
and provide several definitions
that will enable stating the main results.

\subsection{Codimension conjectures}
In general $\mathrm{M}_+(G)$ is difficult to determine exactly,
but some combinatorial bounds are known.
A lower bound on the maximum codimension of a faithful orthogonal
representation was conjectured in 1987 by Maehara.
\newtheorem*{MC}{Maehara's Conjecture}
\begin{MC}[Maehara, credited in \cite{LSS1989}]
\label{cor:maehara}
For any simple graph $G$, a faithful orthogonal representation exists
whose codimension is at least the minimum degree.
Equivalently, the following inequality always holds:
\[
\delta(G) \le \mathrm{M}_+(G)% \ge \delta(G)
.
\]
\end{MC}

\subsubsection{A conjecture from combinatorial matrix theory}
There are conjectures closely related to the conjecture of Maehara that
arose independently (but somewhat later) in the combinatorial
matrix theory community.

In October of 2006 the American Institute of Mathematics (AIM)
hosted a workshop on combinatorial matrix theory \cite{AIM2006}, studying
among other things the maximum nullity problem.
At this workshop a conjecture arose
that, over the reals or other infinite fields,
the minimum degree of a graph should be a lower bound on its maximum nullity.
This inequality
has come to be known as the
Delta Conjecture.\footnote{Despite
the capitalization of the word ``Delta'' in its name,
the conjecture concerns minimum degree~$\delta(G)$
and not maximum degree~$\Delta(G)$.
It is sometimes written as the ``$\delta$-conjecture''.}
Partial results include
the fact that the Delta Conjecture is satisfied for all bipartite graphs \cite{BFHRS2008}
or complements of bipartite graphs \cite{BFHHSv2009}.

\newtheorem*{DC}{Delta Conjecture}
\begin{DC}[Barioli et al.\cite{AIM2006}]%[Friedland \cite{AIM2006}]
For any simple graph $G$, the following inequality holds:
\[
\delta(G) \le \mathrm{M}(G).
\]
\end{DC}

Maehara's Conjecture, which supposes that an orthogonal representation exists with
positive definite Gram matrix, implies the Delta Conjecture.

\subsubsection{A conjecture involving the Strong Arnold Property}

The SAP
was introduced to combinatorial matrix theory
by Colin~de~Verdi\`ere in 1990 in defining the parameter $\mu(G)$ \cite{CdV1990},
which has a number of appealing properties such as the fact that
a graph is planar if and only if $\mu(G) \le 4$.

The SAP fundamentally concerns the
transversality of certain matrix spaces,
but is usually stated in terms of
an equivalent algebraic statement.
The entrywise (or Hadamard, or Schur) product
is written $A \circ B$. %, and $I$ is the identity matrix.

\begin{definition}
\label{def:SAP}
A real symmetric $n \times n$ matrix $A$ is said to have
\emph{the Strong Arnold Property} (\emph{SAP})
if the only real symmetric matrix $X$ that satisfies the three equations
\[
A \circ X = 0,\ \ I \circ X = 0,\ \ \mbox{and }AX = 0
\]
is the matrix $X = 0$.
\end{definition}

Restricting to matrices with SAP
can have useful consequences for nullity parameters.
The parameter $\nu(G)$ is one such, here
defined in terms of orthogonal representations.
\begin{definition}
Given a graph $G$, $\nu(G)$ is the
greatest possible codimension of a faithful orthogonal representation $\varphi$
of $G$ such that the Gram matrix of $\varphi$ has SAP.
\end{definition}

The graph parameter $\nu(G)$ is, like $\mu(G)$, \emph{minor monotone}, which
means that $\nu(F) \le \nu(G)$ holds whenever $F$ is
a contraction of a subgraph of $G$
\cite{CdV1998}\footnote{In the original
proof of the minor monotonicity of $\nu(G)$, the parameter
is called $a(G)$, the ``algebraic width'' of the graph.}.
The inequality $\nu(G) \le \mathrm{M}_+(G)$ suggests
a possible conjecture stronger than that of Maehara.

\newtheorem*{SDC}{Strong Delta Conjecture}
\begin{SDC}[Barioli et al. \cite{BBFHHSvv2012}]
For any simple graph $G$, the following inequality holds:
\[
\delta(G) \le \nu(G).
\]
\end{SDC}
The main results will imply the Strong Delta Conjecture, which
implies Maehara's Conjecture, which implies the Delta Conjecture.

\subsection{The LSS construction}
\label{subsec:LSS}
The earliest literature reference to Maehara's Conjecture is
a study of the codimension problem by
LSS, who
succeeded in showing that under the additional constraint of
general position, $\kappa(G)$ is exactly the maximum codimension.

\begin{definition}
The representation $\varphi$ of $G$ into $\mathbb{R}^\ad$ is in
\emph{general position} if, for any subset $S \subseteq V(G)$ 
with $|S| = k \le \ad$,
the image of $S$ under $\varphi$ is a linearly independent set of $k$ vectors.
\end{definition}

\begin{theorem}[LSS \cite{LSS1989, LSS2000}]
\label{thm:LSS1}
If $G$ is a simple graph with $n$ vertices, then the following are equivalent:
\begin{enumerate}
\item The graph $G$ is $(n-\ad)$-connected.
\item The graph $G$ has a general-position orthogonal representation in $\mathbb{R}^\ad$.
\item \label{generic}The graph $G$ has an orthonormal representation in $\mathbb{R}^\ad$ such that for
each vertex $v$, the vectors representing the vertices nonadjacent to $v$ are linearly independent.
\end{enumerate}
\end{theorem}

Although by this theorem the codimensions in which Condition (2) or
Condition (3) can occur are the same, %namely any codimension less than $\kappa(G)$,
a particular faithful orthogonal representation
may satisfy Condition (3) without satisfying Condition (2).
%as will be shown by example in ...

The LSS construction is probabilistic and depends on the chosen vertex ordering
$V = (v_1, \dots, v_n)$.
For each vertex $v_j$, denote by $W_j$ the order-preserving sublist of $V$ consisting of vertices that both come before and are not adjacent to $v_j$,
whose total number will be denoted $k = k(v_j)$:
\[
W_j = (w_1, w_2, \dots, w_k) := (v_i \in V \ |\  i < j \mbox{ and } v_i \not\sim v_j ).
\]
If $v_j$ comes before all of its neighbors, then $W_j$ will be the empty list $()$ with $k(v_j) = 0$.
The notation $W_j$ will be carried forward, and
%the size of the constraining set $V_j$, to be notated as 
$k(v_j) = |W_j|$ will play a role later on.
The LSS construction chooses representing vectors in order, implying in particular that before
$\varphi(v_j)$ is to be chosen, representations $\varphi(w_1)$ through $\varphi(w_k)$ will already have been chosen.
Denote by $S_j$ the orthogonal complement of their span:
\[
S_j := \{\mathbf{x} \in \mathbb{R}^\ad \ | \ \forall w_i \in W_j,\, \dagg{\mathbf{x}}\varphi(w_i) = 0\}.
\]
The LSS construction simply chooses $\varphi(v_j)$
randomly\footnote{It could for example be chosen uniformly at random
from
Haar measure on the unit sphere in $S_j$ when $S_j$ is nontrivial.
Other reasonable choices give the
same qualitative behavior.}
from within $S_j$.

Note that the ordering of $V$ imposes an asymmetry in the way that random choices are made:
For $i \not\sim j$ with $i < j$, the subspace $S_j$ from which $\varphi(v_j)$ is chosen
depends on $\varphi(v_i)$, but the earlier subspace $S_i$ does not depend on the later choice $\varphi(v_j)$.

\Cref{thm:LSS1} guarantees, in certain dimensions, the existence of an orthogonal representation,
but not necessarily the existence of a faithful orthogonal representation.
A companion theorem in the same work does however guarantee this.

\begin{theorem}[LSS \cite{LSS1989, LSS2000}]
\label{thm:LSS2}
If $G$ is a simple graph with $n$ vertices and
$G$ is $(n - \ad)$-connected,
then the
LSS
construction in $\mathbb{R}^\ad$ yields,
for any ordering of the vertices and with
probability $1$, a faithful orthogonal representation.
\end{theorem}

\begin{corollary}For any simple graph $G$, $\kappa(G) \le \mathrm{M}_+(G)$ holds.
\end{corollary}

\begin{corollary}For any simple graph $G$, $\kappa(G)$ is
the maximum codimension of a general-position faithful
orthogonal representation of $G$.
\end{corollary}

Although the LSS construction itself is straightforward,
the resulting probability
distribution on ensembles of $n$ vectors
exhibits some subtlety.
For example,
the original arguments of \cite{LSS1989}
assumed incorrectly that
the distribution of vector ensembles in the required space is independent
of the choice of vertex ordering.
What is true, and sufficient \cite{LSS2000},
is that distributions under different orderings have
the same sets of measure zero
when the codimension of the representation is at most $\kappa(G)$.

Maehara's Conjecture might be suggested by counting
degrees of freedom.
In dimension $\ad = n - \delta(G)$,
the subspace $S_j$ always has dimension at
least \[\ad - \Delta(\overline{G})  = \ad - (n - \delta(G) - 1) = \ad - \ad + 1 = 1,\]
and so for example the LSS construction in dimension
$n - \delta(G)$ will with probability $1$ yield a
representation by nonzero vectors that can be scaled to an orthonormal representation.
The question is whether the resulting representation
will be faithful, or in other words whether additional,
unwanted orthogonalities can be avoided.
The next example illustrates that ensuring non-trivial subspaces $S_j$,
while it suffices for the existence of an orthonormal representation,
does not suffice for the existence of a faithful representation.

\begin{example}
\label{ex:P4}
Consider the ordered path on four vertices $P_4 = (v_1, v_2, v_3, v_4)$, and
suppose that $G$ is the complement of $P_4$, also a path on four vertices but with different adjacencies
$v_1 \sim v_3$, $v_1 \sim v_4$, and $v_2 \sim v_4$.
Perform LSS in the Cartesian plane $\mathbb{R}^2$, attempting to obtain an orthonormal representation $\varphi$.
The vertex list $W_1$ is empty, so $S_1$ is all of $\mathbb{R}^2$.
Without loss of generality (up to orthogonal rotation), choose
\[\varphi(v_1) = \begin{bmatrix}1 & 0 \end{bmatrix}^\top.\]
The non-adjacency $v_1 \not\sim v_2$ (or equivalently, the adjacency $v_1 \sim v_2$ in $\overline{G}$) gives $W_2 = (v_1)$, so
$S_2$ has dimension $1$, imposing the choice, unique up to sign, of \[ \varphi(v_2)= \begin{bmatrix}0 & 1\end{bmatrix}^\top.\]
In a similar way, $W_3 = (v_2)$ gives $S_3$ of dimension $1$ and imposes, up to sign, \[\varphi(v_3) = \begin{bmatrix}1 & 0 \end{bmatrix}^\top,\]
after which $W_4 = (v_3)$ gives $S_4$ of dimension $1$ and imposes, up to sign, \[\varphi(v_4) = \begin{bmatrix}0 & 1 \end{bmatrix}^\top.\]
At every stage there was at most one previous non-adjacent vertex, and so it was always possible to choose a non-zero vector.
But every choice other than the first is completely constrained up to sign, and regardless
of that first choice, $\varphi(v_4)$ ends up orthogonal to $\varphi(v_1)$, whereas
$v_4 \sim v_1$ in $G$.
An orthonormal representation was possible, but an unwanted orthogonality could not be avoided and a
representation in dimension $\ad = 2$ is never faithful---unlike
in dimension $\ad = n - \delta(G) = 3$, where the Laplacian representation is faithful.
As has long been appreciated, a proof of Maehara's conjecture requires an approach more subtle than just dimension counting.
\end{example}

\subsection{Uniform LSS and weak success}
The subspace $S_j$ from which each vector $ \varphi(v_j)$ is chosen in the LSS construction
has dimension at least $\ad - |W_j|$, and dimension exactly $\ad - |W_j|$ in the case
that the vectors representing $W_j$ are independent.
The fact that the dimension of $S_j$ might vary depending on previous random choices
presents an obstacle to the goal of a uniform parametrization of the LSS construction.
This obstacle is overcome, in a construction that will be called \emph{uniform LSS},
by insisting on the choice $\varphi(v_j) = \mathbf{0}$ whenever the vectors
representing $W_j$ are dependent.
In the case that the vectors representing $W_j$ are independent, but there are $\ad$ of them,
their orthogonal complement $S_j$ is trivial, which also requires the choice $\varphi(v_j) = \mathbf{0}$.
In either case, choosing a zero vector is likely to have cascading consequences, since that same zero vector will
make any later set that contains it dependent.
For each given combination of ambient dimension $\ad$, vertex ordering $(v_1, \dots, v_n)$, and index $j$,
the set of conditions that allow all choices $\varphi(v_j)$ to be nonzero under uniform LSS,
namely that $|W_j| < \ad$ always holds with independent vectors representing $W_j$,
occurs either with probability $0$ or probability $1$.
In the case that it holds with probability $1$, implying
$\varphi(v_j) \ne \mathbf{0}$ for all $j$,
uniform LSS with order $(v_1, \dots, v_n)$ is said to
\emph{succeed weakly} in dimension $\ad$.

\subsection{Strong success of uniform LSS}
\label{subsec:strong-success}
In order for uniform LSS to succeed weakly and also produce a faithful representation, a stronger criterion
for success is required.
For each separate choice of $i$ with $i < j$ and $v_i \sim v_j$ (meaning in particular that $v_i \not \in W_j$),
append $v_i$ (possibly out of order) to the list $W_j$ of length $k = k(v_j)$, and denote the expanded list
\[
W_{j;i} = (w_1, \dots, w_k; v_i).
\]
If uniform LSS succeeds weakly, then with probability $1$ the vectors representing $W_j$ will be independent,
but in the case where $\varphi(v_i)$ is a linear combination of those vectors,
$\varphi(v_j)$ chosen from $S_j$ will be orthogonal to $\varphi(v_i)$ even though
$v_i$ and $v_j$ are adjacent in $G$, which means that
the representation will not be faithful.
Supposing that uniform LSS does succeed weakly,
for each $i < j$ with $v_i \sim v_j$
the vectors representing $W_{j;i}$ will be independent either with probability $0$ or probability $1$.
In the case that uniform LSS succeeds weakly and that in addition
each $W_{j;i}$ produces independent vectors with probability $1$,
uniform LSS is said to \emph{succeed strongly}.

\begin{proposition}
    \label{prop:strongly}
    Given a graph $G$ on $n$ vertices, vertex ordering $(v_1, \dots, v_n)$, and ambient dimension $\ad$,
    uniform LSS succeeds strongly if and only if
    \begin{itemize}
        \item uniform LSS succeeds weakly and 
        \item uniform LSS produces, with probability $1$, a faithful orthogonal representation.
    \end{itemize}
\end{proposition}
\begin{proof}
    Suppose that uniform LSS succeeds weakly and does, with probability $1$, for each $j$ and each $i < j$ such that $v_i \sim v_j$ and $v_i \not \in W_j$, find $\varphi(v_j)$ that is not orthogonal to $\varphi(v_i)$.
    It follows directly that $\varphi(v_i)$ is not in the span of the independent vectors representing $W_j$, each of which is orthogonal to $\varphi(v_j)$, and hence that the vectors representing $W_{j;i}$ are also independent.
    
    Conversely, if uniform LSS succeeds strongly then by definition it also succeeds weakly, which takes care of the first clause.
    Call $S_{j;i}$ the orthogonal complement of the span of the vectors in $W_{j;i}$.
    The subspace $S_{j;i}$ has dimension $\ad - |W_{j;i}|$ with probability $1$ in the case that uniform LSS succeeds strongly.
    As before, call $S_j$ the orthogonal complement of the span of the vectors in $W_j$.
    The subspace $S_j$ has dimension $\ad - |W_j|$ in the case that uniform LSS succeeds at least weakly.
    Given
    \[
        \ad - |W_j| = \ad - (|W_{j;i}| - 1) = (\ad - |W_{j;i}|) + 1,
    \]
    it follows that when $\varphi(v_j)$ is chosen randomly in $S_j$ there is probability $0$ that it is chosen
    from the codimension-$1$ subspace $S_{j;i}$,
    and therefore that with probability $1$ the representation chosen by uniform LSS is faithful.
\end{proof}

\Cref{ex:P4} above provides an illustration of uniform LSS that succeeds weakly but does not succeed strongly:
When choosing the last representing vector $\varphi(v_4)$,
the vector set representing $W_4 = (v_3)$ is independent with probability $1$,
but for the adjacent pair of vertices $v_1 \sim v_4$,
the vector set representing $W_{4;1} = (v_3; v_1)$ is dependent with probability $1$,
and correspondingly $\varphi(v_4)$ is chosen, unfaithfully, orthogonal to $\varphi(v_1)$.

\subsection{Weaker than general position but stronger than SAP}
When uniform LSS succeeds strongly, this will imply, with probability $1$, not only that the
representation is faithful, but that the resulting Gram matrix satisfies
a certain genericity condition that in turn implies SAP.

\begin{definition}
\label{def:upperzero}
A real symmetric matrix $A$ is said to be
\emph{upper-zero generic} if the
diagonal entries of $A$ are nonzero and
if, in each column of $A$,
the zeros above the diagonal
determine an independent set of rows of $A$.
\end{definition}

\begin{proposition}[Strong success implies upper-zero generic]
\label{prop:upperzero}
Given a graph $G$ with vertex ordering $(v_1, \dots, v_n)$,
if uniform LSS in ambient dimension $\ad$ succeeds strongly, then with probability $1$ the
resulting representation is faithful with a Gram matrix $A$ that is upper-zero generic.
\end{proposition}

\begin{proof}
Throughout this proof, whenever a condition is stated to hold,
what is actually meant is that it holds with probability $1$.

The fact that uniform LSS succeeds at least weakly implies that $\varphi(\sigma_j)$ is
not~$\mathbf{0}$, and therefore that none of the diagonal entries 
$a_{jj} = \dagg{\varphi(v_j)} \varphi(v_j)$
of $A$ is~$0$.
Faithfulness when uniform LSS succeeds strongly is by \cref{prop:strongly}, and implies
that the zeros above the diagonal entry in column $j$ of $A$ 
correspond exactly to those indices $i$ such that $v_i \in W_j$.
Since the vectors representing $W_j$ are independent,
the submatrix $A[W_j]$ is invertible,
and those rows determine an independent set of rows of $A$.
Since $j$ was arbitrary,
in every column the set of zeros above the diagonal determines an independent set of rows of $A$, 
implying that $A$ is upper-zero generic.
\end{proof}

\begin{theorem}[Upper-zero generic implies SAP]
\label{thm:upper_to_SAP}
Given a real symmetric matrix $A$ that is upper-zero generic, $A$ satisfies the Strong Arnold Property.
\end{theorem}

\begin{proof}
Suppose that real symmetric $A$ of size $n \times n$ is upper-zero generic,
and let $X$ be an $n \times n$ real symmetric matrix satisfying $A \circ X = 0$ (which implies $I \circ X = 0$) and $AX = 0$.
Suppose for the sake of contradiction that $X$ is not the zero matrix, and
let $j$ be the largest index of a non-zero column $\mathbf{x}_j$ of $X$.
By the symmetry of $X$ and the fact $I \circ X = 0$, the support of $\mathbf{x}_j$ lies
strictly above the diagonal of $X$, and by the fact $A \circ X = 0$, the support of $\mathbf{x}_j$
lies in the positions of zero entries above the diagonal in column $j$ of $A$.
Since $A$ is upper-zero generic, this set of indices corresponds to an independent set of rows
of $A$, and also, by the symmetry of $A$, corresponds to an independent set of columns of $A$.
Thus $A\mathbf{x}_j$, which is column $j$ of $AX=0$, is both a zero vector and
a nonzero linear combination of an independent set of columns of $A$.
The conclusion, by contradiction, is that $X$ is in fact the zero matrix, which demonstrates that $A$ satisfies SAP.
\end{proof}

\begin{definition}
\label{def:Mseq}
Let $G$ be a graph on $n$ vertices. The  \emph{maximum PSD upper nullity} of $G$,
denoted $\Mpu(G)$, is the largest nullity among all matrices $A = [a_{ij}]$
such that
\begin{enumerate}
\item the matrix $A$ of size $n \times n$ is real symmetric and positive semidefinite;
\item there exists a permutation matrix $P$ such that $P^{-1}AP$ is upper-zero generic;
\item for 
$i \ne j$, $a_{ij} = 0$ holds whenever $v_i \not\sim v_j$ in $G$; and
\item for
$i \ne j$, $a_{ij} \ne 0$ holds whenever $v_i \sim v_j$ in $G$.
\end{enumerate}
\end{definition}

\begin{proposition}
\label{prop:seq}
For any simple graph $G$, $\Mpu(G) \le \nu(G)$ holds.
\end{proposition}

\begin{proof}
Let $A$ be a matrix of nullity $\Mpu(G)$
satisfying the hypotheses of \cref{def:Mseq}.
By \cref{thm:upper_to_SAP}, the matrix $P^{-1}AP$ satisfies SAP,
and since SAP is permutation invariant, so does $A$,
making $A$ itself a matrix of nullity sufficient to show that $\Mpu(G) \le \nu(G)$.
\end{proof}

\subsection{Greedegree}
\label{subsec:greedegree}
The main results will show that faithful orthogonal representations
with Gram matrices that are upper-zero generic exist, not
only in codimension equal to $\delta(G)$, but also in
codimension
equal to a new degree parameter,
to be called the greedegree of $G$, that can
in general be larger than $\delta(G)$.

A maximum cardinality search (MCS) \cite{TY1984} visits every vertex of a graph once,
with the constraint that each new vertex chosen must have a maximum cardinality of adjacencies to previously visited vertices.
The process can equivalently be viewed as one in which the graph is built up, vertex by vertex, as a nested series of induced subgraphs; the MCS constraint is then that each new vertex must greedily maximize the number of new edges added to the graph at each stage.
Adopting this viewpoint, an MCS ordering of a graph will herein usually be referred to as a \emph{greedy} ordering.
A rigorous statement of greedy ordering is facilitated by terminology that covers both the number of edges that are actually added at stage $i$,
and every alternate number of edges that could have been added by choosing vertex $v_j$ instead of $v_i$.

\begin{definition}
\label{def:local}
Given an ordering
$V(G) =(v_1, \dots, v_{n})$
and integers $i$ and $j$ such that $1 \le i \le j \le n$,
the \emph{local degree of $v_j$ at stage $i$} is the number of
edges between $G\left[\left\{v_1, \dots, v_{i - 1}\right\}\right]$ and $v_j$,
or in other words the size of the set $\{v_hv_j \in E \st h < i\}$.
\end{definition}

\begin{definition}
\label{def:greedy}
An ordering $V(G) =(v_1, \dots, v_{n})$ is called MCS or \emph{greedy}
if, for all $1 \le i \le j \le n$,
the local degree of $v_i$ at stage $i$ is at least as great as the
local degree of $v_j$ at stage $i$.
\end{definition}

In common with many other algorithms termed ``greedy'', the maximizing constraint is only local to each step,
and two greedy orderings of the same graph may differ in the total number of induced edges at
intermediate stages.
The local constraint does however result in global tendencies to produce dense early stages,
to postpone low degree vertices until later in the process,
and to choose an especially low-degree vertex as the final vertex of the ordering.
Of particular interest is how well these global tendencies can be avoided,
as measured by having a final vertex of relatively high degree.
For example, for $n \ge 2$ any single low degree vertex can be avoided as the final vertex simply by choosing it as the first vertex.
Similarly, for $n \ge 3$ any adjacent pair of low degree vertices can both be avoided as the final vertex simply by visiting them first and second in a maximum cardinality search.
\begin{definition}
For $G$ a simple graph, the
\emph{maximal greedy final degree}
of $G$
or \emph{greedegree}
of $G$,
denoted $\gdel(G)$,
is the maximum, over all greedy orderings of $G$, of the degree of the final vertex in the ordering.
\end{definition}

\begin{remark}
It is known \cite{Lucena2003} that in any maximum cardinality search of $G$, the maximum that occurs over all stages of the number of previously visited adjacent vertices is a lower bound on the treewidth of $G$, denoted $\mathrm{tw}(G)$.
From this it follows that the maximum number of edges that can added at the end of a greedy ordering also gives a lower bound,
\[
\gdel(G) \le \mathrm{tw}(G).
\]
\end{remark}

\begin{example}
    \label{ex:tree-gdg}
    Let $T$ be a tree on $n \ge 2$ vertices.
    The lowest degree vertices of $T$ are vertices of degree $1$, and unavoidably one of these must end up as the final vertex of any greedy ordering, because at each intermediate stage, and the penultimate stage in particular, the induced graph is always connected.
    It follows that $\gdel(T) = 1$.
\end{example}

The obvious bounds are $\delta(G) \le \gdel(G) \le \Delta(G)$, with equality in both
cases for a regular graph.
In typical cases strictness of the upper inequality cannot be avoided,
but in some cases (not including trees) it is at least possible to make the lower inequality strict as well.
A general strategy for obtaining large final degree that performs better
than exhaustive search is not
known, and the calculation of $\gdel(G)$ is shown in \cref{sec:NP} to be an NP-hard problem---in fact, it is further shown to be NP-hard to determine, for an arbitrary pair of vertices that are {\em not} adjacent,
whether they can both be avoided as the final vertex of the ordering.

\subsection{Statement of main results}
\label{subsec:mainstate}

Much of the section up to this point can be summarized by two chains of inequalities
\[ \kappa(G) \le \delta(G) \le \gdel(G) \]
and
\[ \Mpu(G) \le \nu(G) \le \M_{+}(G) \le \M(G). \]
The Main Theorem, to be proven in \cref{subsec:mainproof}, will allow these to be joined into a single chain of inequalities.

\begin{restatable}[Main Theorem]{theorem}{mainthm}
\label{thm:main}
Let $G$ be a simple graph whose $n$ vertices
$V(G) = (v_1, \dots, v_n)$
are ordered according to a greedy ordering,
and let $\ad = n - \mathrm{deg}(v_n)$.
Then uniform LSS
on $(v_1, \dots, v_n)$ in $\mathbb{R}^\ad$ succeeds strongly.
\end{restatable}

Combining the Main Theorem with \cref{prop:upperzero} and \cref{thm:upper_to_SAP}, the following corollaries are implied.

\begin{corollary}
For any simple graph $G$, $\gdel(G) \le \Mpu(G)$.
\end{corollary}

\begin{corollary}
The Strong Delta Conjecture holds.
\end{corollary}

\begin{corollary}
Maehara's Conjecture holds.
\end{corollary}

\begin{corollary}
The Delta Conjecture holds.
\end{corollary}

\begin{remark}
    The parameter $\gdel(G)$ is not minor monotone, but SAP-based nullity parameters such as $\nu(G)$ are.
    This can improve the lower bound:
    If $G^\prime$ is a minor of $G$ that breaks monotonicity by exhibiting $\gdel(G^\prime) > \gdel(G)$, then
    \[\gdel(G^\prime) \le \nu(G^\prime) \le \nu(G)\] gives a strictly stronger bound than
    $\gdel(G) \le \nu(G)$.
\end{remark}

\section{The hanging garden construction}
\label{sec:hanging}
\subsection{Parametrization}
\label{subsec:parametrization}
The goal is to completely parametrize all of the random choices that are made during uniform LSS using free variables that encode the choices made at each stage.
At the end of the process, each entry $r_{ij}$ of the matrix $R$---which is to say, each entry $i$ of every vector $\varphi(v_j)$---and each entry $t_{ij}$ of the overall Gram matrix $T$ will be a (large, complicated) polynomial with integer coefficients over a set of independent commuting variables.
The hanging garden construction will organize these large polynomials diagrammatically in a way that will facilitate reasoning about the main results.

\subsection{Ordering and dependencies}
It has been assumed that $G$ has
a fixed ordering 
\[
V(G) = (v_1, \dots, v_n).
\]
It will not yet be assumed that the ordering is greedy (although that assumption will come later).
Adjacent vertices of $G$ do not directly affect each other in the LSS construction,
while nonadjacent vertices $i \not \sim j$ have an asymmetric relationship that depends on the ordering:
For $i < j$, with $v_i \in W_j$ but $v_j \not \in W_i$,
the earlier chosen $\varphi(v_i)$ imposes an orthogonality restraint on the later choice of $\varphi(v_j)$, but not the other way around.
The number of orthogonality constraints that uniform LSS imposes on $\varphi(v_j)$ is denoted by $k$,
recalling that
\[
k(v_j) = |W_j|.
\]

\subsection{Independent and dependent node variables}
\label{subsec:nodevar}
Each vertex $v_1, v_2, v_3, \dots$ of $G$ gives rise to a class $\alpha, \beta, \gamma, \dots$ of node variables,
each with an associated $k = k(v_j)$ that determines an \emph{arity} (number of arguments) equal to $k + 1$.
A \emph{class of node variables} is an alternating function 
\[\{1, \dots, \ad\}^{k+1} \rightarrow \mathbb{R}.\]
The arguments to the function are indices, and the resulting variable is written using subscript notation rather than functional notation.
To give an example, suppose that $k(v_3) = 2$ gives rise to class $\gamma$ of node variables of arity $2 + 1$.
Instead of writing $\gamma(s, t, u)$ when applying three index arguments to the alternating function $\gamma$, the result is written
\[\gamma_{s, t; u}.\]
(Some explanation is perhaps due: Why express arity $3$ as the sum $2 + 1$, and why place a semicolon before the last subscript of $\gamma$? These choices foreshadow that, inside a garden diagram, the last argument of an alternating function such as $\gamma$ is special.)
The \emph{independent node variables} are those for which the indices are in strictly increasing order.
The number of free variables or independent degrees of freedom within a class of node variables is thus $\binom{\ad}{k + 1}$.
For example, $\ad = 4$ would give the following independent node variables within the above class $\gamma$:
\[
\gamma_{1, 2; 3}, \gamma_{1, 2; 4}, \gamma_{1, 3; 4}, \mbox{ and } \gamma_{2, 3; 4}.
\]
The function that determines a class of node variables must be alternating in a sense that is equivalent to these two constraints:

(1) Any node variable in which an index is repeated evaluates to $0$; for example, \[\gamma_{3, 1; 3} := 0.\]

(2) Any node variable in which the indices are distinct, but not in increasing order, is a \emph{dependent node variable} and evaluates either directly to an independent node variable, or to its negation, with a sign given by the sign of the permutation that puts the arguments in order.
For example,
\begin{eqnarray*}
\gamma_{2, 1; 3} &:=& -\gamma_{1, 2; 3}\\
\mbox{ and }\ \ \gamma_{4, 1; 3} &:=& \gamma_{1, 3; 4}
\end{eqnarray*}
are two dependent node variables.

\subsection{Gardens of Type I and Type II}
\begin{figure}
\begin{tikzpicture}[
  font=\normalsize,
  conduit/.style={line width=1.2pt, line cap=round, line join=round},
  garden/.style={draw, fill=black!25},
  input/.style={draw, fill=white, minimum height=0.60cm, inner sep=1pt},
  nodebox/.style={draw, rounded corners=8pt, fill=black!10,
                  minimum height=0.70cm, inner sep=0pt}
]
  % ------------------------------------------------------------
  % 1) Garden boundary/background (extra top headroom)
  % ------------------------------------------------------------
  \draw[garden] (0.20,0.35) rectangle (7.40,8.50);

  % Column x-positions (tweaked)
  \def\xIone{1.3}   % input 1 column
  \def\xA{1.35}      % beta1/alpha column
  \def\xB{2.30}      % theta column
  \def\xC{3.00}      % gamma column moved further left (was 3.45)
  \def\xD{5.25}      % right delta / alpha column
  \def\xE{6.20}      % far-right input/beta column
  \def\xAtop{3.55}   % top alpha moved left (was 3.90)

  % Row y-positions
  \def\yOne{7.25}
  \def\yTwo{6.05}
  \def\yThree{4.80}
  \def\yFour{3.55}
  \def\yFive{2.35}
  \def\ySix{1.15}

  % ------------------------------------------------------------
  % Component centers
  % ------------------------------------------------------------
  % Row 1
  \coordinate (TH) at (\xB,\yOne);
  \coordinate (ZE) at (5.60,\yOne);

  % Row 2
  \coordinate (I1)  at (\xIone,\yTwo);
  \coordinate (Aup) at (\xAtop,\yTwo);
  \coordinate (I4)  at (\xE,\yTwo);

  % Row 3 (deltas)
  \coordinate (D1) at (2.10,\yThree);
  \coordinate (D2) at (5.55,\yThree);

  % Row 4
  \coordinate (GA)  at (\xC,\yFour);   % gamma now under left delta
  \coordinate (Amd) at (\xD-0.3,\yFour);

  % Row 5
  \coordinate (B1) at (\xA,\yFive+0.3);
  \coordinate (I2) at (\xC+0.1,\yFive);
  \coordinate (B2) at (\xE,\yFive);

  % Row 6
  \coordinate (Abt) at (\xA,\ySix+0.3);
  \coordinate (I3)  at (\xE,\ySix);

  % ------------------------------------------------------------
  % 2) Conduits
  % ------------------------------------------------------------

  % Type II cap conduit between theta and zeta (inside garden)
  \draw[conduit, rounded corners=10pt]
    ($(TH)+(0.2,0.35)$) -- ($(TH)+(0.2,0.85)$)
    -- ($(ZE)+(-0.2,0.85)$) -- ($(ZE)+(-0.2,0.00)$);

  % theta lower conduits (3)
  \draw[conduit] (\xIone,\yOne-0.35) -- (\xIone,\yTwo+0.00);  % to I1
  \draw[conduit] (\xB,  \yOne-0.35) -- (\xB,  \yThree+0.00); % to D1 (vertical above delta)
  \draw[conduit] (\xAtop-0.2,\yOne-0.35) -- (\xAtop-0.2,\yTwo+0.00); % to Aup

  % zeta lower conduits (2)
  \draw[conduit] (5.05,\yOne-0.35) -- (5.05,\yThree+0.00);   % to D2 (vertical above delta)
  \draw[conduit] (\xE,\yOne-0.35) -- (\xE,\yTwo+0.00);        % to I4

  % Delta 1 lower conduits (2)
  \draw[conduit] (\xA,\yThree-0.35) -- (\xA,\yFive+0.00);     % lower 1 -> B1 (passthrough row 4)
  \draw[conduit] (\xC-0.1,\yThree-0.35) -- (\xC-0.1,\yFour+0.00);     % lower 2 -> gamma

  % Gamma lower conduit
  \draw[conduit] (\xC+0.1,\yFour-0.35) -- (\xC+0.1,\yFive+0.00);      % -> I2

  % Beta1 lower conduit
  \draw[conduit] (\xA,\yFive-0.05) -- (\xA,\ySix+0.30);       % -> Abt

  % Delta 2 lower conduits (2)
  \draw[conduit] (\xD-0.4,\yThree-0.35) -- (\xD-0.4,\yFour+0.00);     % lower 1 -> Amd
  \draw[conduit] (\xE,\yThree-0.35) -- (\xE,\yFive+0.00);     % lower 2 -> B2 (passthrough row 4)

  % Beta2 lower conduit
  \draw[conduit] (\xE,\yFive-0.35) -- (\xE,\ySix+0.00);       % -> I3

  % ------------------------------------------------------------
  % 3) Components (draw last)
  % ------------------------------------------------------------

  % Row 1 nodes
  \node[nodebox, minimum width=3.10cm] at (TH) {$\theta$};
  \node[nodebox, minimum width=2.60cm] at (ZE) {$\zeta$};

  % Row 2
  \node[input,  minimum width=0.95cm] at (I1)  {$\mathcal{I}_1$};
  \node[nodebox, minimum width=1.50cm] at (Aup) {$\alpha$};
  \node[input,  minimum width=0.95cm] at (I4)  {$\mathcal{I}_4$};

  % Row 3 deltas
  \node[nodebox, minimum width=2.80cm] at (D1) {$\delta$};
  \node[nodebox, minimum width=2.80cm] at (D2) {$\delta$};

  % Row 4
  \node[nodebox, minimum width=1.50cm] at (GA)  {$\gamma$};
  \node[nodebox, minimum width=1.50cm] at (Amd) {$\alpha$};

  % Row 5
  \node[nodebox, minimum width=1.50cm] at (B1) {$\beta$};
  \node[input,  minimum width=0.95cm] at (I2) {$\mathcal{I}_2$};
  \node[nodebox, minimum width=1.50cm] at (B2) {$\beta$};

  % Row 6
  \node[nodebox, minimum width=1.50cm] at (Abt) {$\alpha$};
  \node[input,  minimum width=0.95cm] at (I3)  {$\mathcal{I}_3$};

\end{tikzpicture}

\caption{A hanging garden diagram of Type II}
\label{fig:garden-ii}
\end{figure}
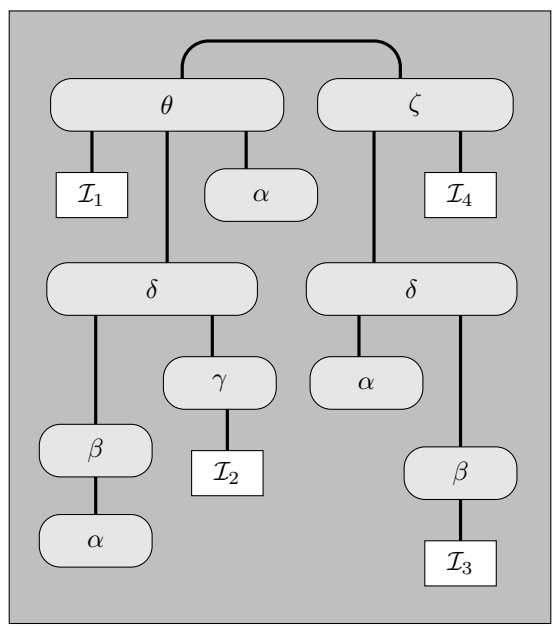

A \emph{hanging garden diagram} (or just \emph{garden,} for short) takes the form of a rectangle with various internal components connected by conduits, as illustrated by \cref{fig:garden-i} (a hanging garden of Type~I) and \cref{fig:garden-ii} (a hanging garden of Type~II).
Every conduit connects below to the upper border of at least one internal component, and almost every conduit forms a vertical connection that also connects above to the lower border of another internal component.
The single exception to this is one uppermost conduit within the garden.
The uppermost conduit connects as follows:
\begin{itemize}
    \item A garden of Type I has an uppermost conduit that exits the upper border of the rectangular garden as a single external connection.
    \item A garden of Type II has an uppermost conduit in the form of a $\cap$-connection between the upper borders of two internal components. Type II gardens have no external connections.
\end{itemize}

\subsection{Internal components and conduit placement}
Every internal component of a garden takes one of two possible forms:
\begin{itemize}
    \item an \emph{input box} is a rectangle-shaped void within the diagram, with a single conduit that exits vertically from its upper border.
    The inputs are numbered $\mathcal{I}_1$ through $\mathcal{I}_q$, left to right---not according to their physical position, but lexicographically, with higher branching taking precedence for the ordering.
    \item a \emph{node} is a rounded rectangle that is labeled by the name of some node variable class (e.g., $\alpha, \beta, \gamma, \dots$).
    Every node has a single conduit that exits vertically from its upper border.
    In addition, for a node variable class of arity $k + 1$, the node has $k$ conduits that descend vertically from its lower border, numbered from left to right as conduits $1$ through $k$.
    The upper conduit is numbered as conduit $k + 1$.
    Multiple nodes within the same garden may be labeled by the same class of node variables, and nodes that share a name will also share the same arity.

\end{itemize}

\begin{figure}
\begin{tikzpicture}[
  font=\normalsize,
  conduit/.style={line width=1.2pt, line cap=round, line join=round},
  garden/.style={draw, fill=black!25},
  input/.style={draw, fill=white, minimum height=0.60cm, inner sep=1pt},
  nodebox/.style={draw, rounded corners=8pt, fill=black!10,
                  minimum height=0.70cm, inner sep=0pt}
]

  % ------------------------------------------------------------
  % 1) Garden boundary/background (Type I: top conduit exits)
  % ------------------------------------------------------------
  \draw[garden] (0.00,0.20) rectangle (6.9,7.8);

  % ------------------------------------------------------------
  % Component centers
  % ------------------------------------------------------------
  \coordinate (TH) at (3.40,7.10);  % theta (top, wide)

  \coordinate (I3) at (3.60,6.05);  % I3 (row below)

  % Reverted (preferred) placement/width for zeta and upper beta:
  \coordinate (ZE) at (2.50,5.10);  % zeta (wide-ish, left)
  \coordinate (BE) at (5.60,5.10);  % beta (upper, narrow, right)

  \coordinate (I1) at (1.30,4.05);  % I1 (below zeta, left)
  \coordinate (DE) at (3.60,4.05);  % delta (slightly less wide than before)

  \coordinate (GA) at (2.40,2.95);  % gamma
  \coordinate (AR) at (5.60,2.95);  % alpha (right, under upper beta)

  \coordinate (AL) at (2.40,1.75);  % alpha (left, under gamma)
  \coordinate (B2) at (4.20,2.05);  % beta (lower, under delta's right conduit)

  \coordinate (I2) at (4.20,1.00);  % I2 (bottom)

  % Attachment x-positions for multi-conduit nodes
  \coordinate (thL) at (2.50,0);    % theta -> zeta (aligned with ZE)
  \coordinate (thM) at (3.40,0);    % theta -> I3
  \coordinate (thR) at (5.60,0);    % theta -> upper beta

  \coordinate (zeL) at (1.30,0);    % zeta -> I1 (left lower conduit)
  \coordinate (zeR) at (3.80,0);    % zeta -> delta (right lower conduit, chosen to the right)

  \coordinate (deL) at (2.80,0);    % delta -> gamma
  \coordinate (deR) at (3.80,0);    % delta -> lower beta

  % ------------------------------------------------------------
  % 2) Conduits (slightly over-long; components drawn later hide overlap)
  % ------------------------------------------------------------

  % Type I unique uppermost conduit exiting the garden
  \draw[conduit] ($(TH)+(0,0.35)$) -- (3.40,8.10);

  % Theta's three lower conduits
  \draw[conduit] ($(thM)+(0.2,6.75)$) -- ($(I3)+(0,0.00)$);  % theta -> I3
  \draw[conduit] ($(thL)+(-1,6.75)$) -- ($(ZE)+(-1,0.00)$);  % theta -> zeta
  \draw[conduit] ($(thR)+(0,6.75)$) -- ($(BE)+(0,0.00)$);  % theta -> beta (upper)

  % Zeta's two lower conduits (vertical to their targets via chosen x-attachments)
  \draw[conduit] ($(zeL)+(0,4.75)$) -- ($(I1)+(0,-0.05)$); % zeta -> I1
  \draw[conduit] ($(zeR)+(0,4.75)$) -- ($(DE)+(0.2,0.00)$);  % zeta -> delta

  % Upper beta down to right alpha
  \draw[conduit] ($(BE)+(0,-0.35)$) -- ($(AR)+(0,0.00)$);

  % Delta down to gamma (left) and to lower beta (right)
  \draw[conduit] ($(deL)+(0,3.70)$) -- ($(GA)+(0.4,0.00)$);
  \draw[conduit] ($(deR)+(0.4,3.70)$) -- ($(B2)+(0,0.00)$);

  % Gamma down to left alpha
  \draw[conduit] ($(GA)+(0,-0.35)$) -- ($(AL)+(0,0.00)$);

  % Lower beta down to I2
  \draw[conduit] ($(B2)+(0,-0.35)$) -- ($(I2)+(0,-0.05)$);

  % ------------------------------------------------------------
  % 3) Components (draw last)
  % ------------------------------------------------------------

  % theta (wide, spanning garden comfortably)
  \node[nodebox, minimum width=6.00cm] at (TH) {$\theta$};

  % inputs
  \node[input, minimum width=1.10cm] at (I3) {$\mathcal{I}_3$};
  \node[input, minimum width=0.95cm] at (I1) {$\mathcal{I}_1$};
  \node[input, minimum width=0.95cm] at (I2) {$\mathcal{I}_2$};

  % middle layer nodes (zeta wide; upper beta narrow)
  \node[nodebox, minimum width=4.20cm] at (ZE) {$\zeta$};
  \node[nodebox, minimum width=1.50cm] at (BE) {$\beta$};

  % delta (slightly less wide)
  \node[nodebox, minimum width=3.00cm] at (DE) {$\delta$};

  % gamma and right alpha
  \node[nodebox, minimum width=2.20cm] at (GA) {$\gamma$};
  \node[nodebox, minimum width=1.50cm] at (AR) {$\alpha$};

  % left alpha and lower beta (both match upper beta width)
  \node[nodebox, minimum width=1.50cm] at (AL) {$\alpha$};
  \node[nodebox, minimum width=1.50cm] at (B2) {$\beta$};

\end{tikzpicture}

\caption{A hanging garden diagram of Type I}
\label{fig:garden-i}
    
\end{figure}
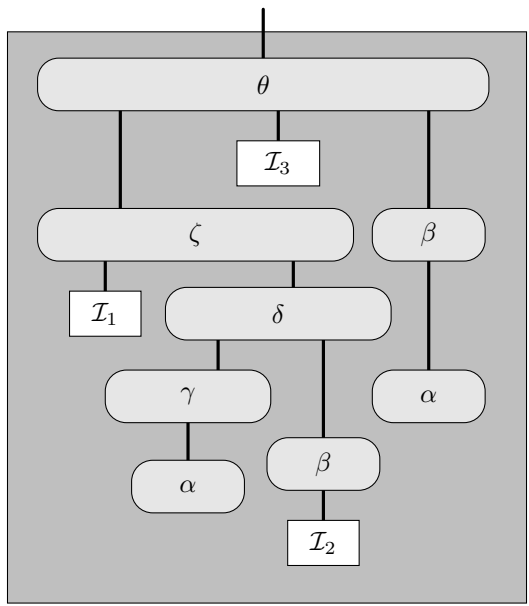

\subsection{Input matching, conduit indexing, and valuation}
The vertical placement of the $q$ internal input boxes of a garden does not matter, but only their ordering $\mathcal{I}_1, \dots, \mathcal{I}_q$.
A garden is \emph{matched} if each input $\mathcal{I}_j$
has been provided with an input vector $\mathbf{x}_j \in \mathbb{R}^\ad$.
(A garden with no internal input boxes, $q=0$, is vacuously matched without any choice of vectors.)
A matched garden is \emph{indexed} if every conduit, including the uppermost, has been labeled by an index from the set $\{1, \dots, \ad\}$.
The \emph{valuation} of an indexed matched garden is calculated as a product of \emph{factors} as follows:
\begin{itemize}
    \item Every internal input box $\mathcal{I}_j$ is connected at its upper border to a single conduit with index label $i$, and has been provided with an input vector $\mathbf{x}_j$.
    This input box contributes a factor that is entry $i$ of vector $\mathbf{x}_j$.
    \item Every internal node is labeled by the name of a node variable class, say $\alpha$, of arity $k+1$.
    It has $k$ lower conduits with index labels (in order, left-to-right) $s_1$ through $s_k$, and an upper conduit with index label $s_{k+1}$.
    This internal node contributes a factor that is equal to the node variable evaluation
    \[\alpha_{s_1, \dots, s_k; s_{k+1}},\]
    which is either an independent node variable, a dependent node variable that evaluates to a signed independent node variable, or $0$ if some index is repeated.
    \item An indexed matched garden of Type I has, exiting from its upper border, a single conduit with index label $i$.
    This uppermost connection contributes a \emph{vector factor} that is equal to the $i^{\mathrm{th}}$ standard basis vector $\mathbf{e}_i \in \mathbb{R}^\ad$.
    \item The overall valuation of the labeled matched garden is the product of all contributed factors, which for a garden of Type I is a single-support vector in $\mathbb{R}^\ad$, and which for a garden of Type II is a scalar.
\end{itemize}
Note that whenever it happens that some internal node has two adjacent conduits (possibly including its upper conduit number $k + 1$) with the same index label, the overall valuation of the indexed matched garden is zero (i.e., the zero vector in Type I, or the scalar $0$ in Type II).

\begin{observation}
    \label{obs:multi-val}
    The valuation of an indexed matched garden of Type~I or of Type~II is multilinear in each input vector of the matching.
\end{observation}

\subsection{Gardens as vector functions}

Every garden $g$ of Type~I with $q$ input boxes can operate as a function from vector $q$-tuples to vectors,
\[
g: \left( \mathbb{R}^\ad \right)^q \longrightarrow \mathbb{R}^\ad.
\]
Similarly, every garden $g$ of Type~II with $q$ input boxes can operate as a function from vector $q$-tuples to scalars,
\[
g: \left( \mathbb{R}^\ad \right)^q \longrightarrow \mathbb{R}.
\]
In either case, let $p$ be the total number of conduits in the garden (including the uppermost conduit).
The value of the function on inputs $\mathbf{x}_1, \dots, \mathbf{x}_q$ is obtained by matching the input vectors to the input boxes, in order, and then
summing every possible valuation, over all $\ad^p$ ways that the conduits of the matched garden can be indexed.

\begin{observation}
    \label{obs:multilinear}
    The multilinearity of valuation of an indexed matched garden diagram extends to multilinearity, over all inputs,
    of the function defined by a hanging garden diagram of Type~I or of Type~II.
\end{observation}

\subsection{Zero-input gardens as vectors or scalars}
A garden with $q=0$ input boxes gives a function with a one-point domain, which can be evaluated directly.
A zero-input garden of Type~I represents a vector in $\mathbb{R}^\ad$, each entry of which is an integer polynomial over the independent node variables.
A zero-input garden of Type~II represents a scalar that is an integer polynomial over the independent node variables.

\begin{example}
    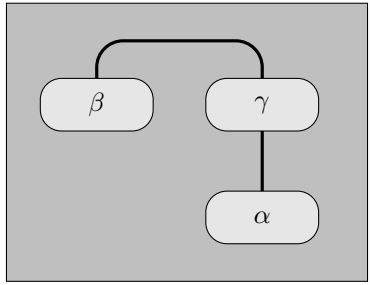
\begin{figure}

\begin{tikzpicture}[
  font=\normalsize,
  conduit/.style={line width=1.2pt, line cap=round, line join=round},
  garden/.style={draw, fill=black!25},
  nodebox/.style={draw, rounded corners=8pt, fill=black!10,
                  minimum height=0.70cm, inner sep=0pt}
]
  % 1) Garden
  \draw[garden] (0.20,1.10) rectangle (5.00,4.80);
  
  % Node centers
  \coordinate (B) at (1.40,3.45); % beta (top-left)
  \coordinate (G) at (3.60,3.45); % gamma (top-right)
  \coordinate (A) at (3.60,1.95); % alpha (below gamma)

  % 2) Conduits
  % Type II cap conduit (beta to gamma)
  \draw[conduit, rounded corners=10pt]
    ($(B)+(0,0.35)$) -- ($(B)+(0,0.85)$)
    -- ($(G)+(0,0.85)$) -- ($(G)+(0,0.00)$);

  % Vertical conduit gamma -> alpha (over-long)
  \draw[conduit] ($(G)+(0,-0.35)$) -- ($(A)+(0,0.00)$);

  % 3) Nodes
  \node[nodebox, minimum width=1.50cm] at (B) {$\beta$};
  \node[nodebox, minimum width=1.50cm] at (G) {$\gamma$};
  \node[nodebox, minimum width=1.50cm] at (A) {$\alpha$};

\end{tikzpicture}
\caption{A zero-input hanging garden of Type~II, which evaluates to an integer polynomial over independent node variables.}
\label{fig:garden-simple-example}
    \end{figure}
Consider the zero-input garden of Type II depicted in \cref{fig:garden-simple-example}.
There are only two conduits.
The uppermost conduit connects on the left to a node labeled $\beta$, of arity $0 + 1$, and on the right to a node labeled $\gamma$, of arity $1 + 1$.
Below $\gamma$ is a node labeled $\alpha$, of arity $0 + 1$.
Take $\ad = 2$.
There are four possible ways to index the conduit pair and get a scalar valuation.
The sum of those valuations,
\[
0 + \alpha_1 \gamma_{2; 1} \beta_2 + \alpha_2 \gamma_{1; 2} \beta_1 + 0 =
(\alpha_2\beta_1 - \alpha_1\beta_2)\gamma_{1; 2},
\]
is, as promised, an integer polynomial over the independent node variables.
\end{example}

\begin{observation}
   \label{obs:multihom}
    The scalar polynomial of a zero-input garden of Type~II and the vector entry polynomials of a zero-input garden of Type~I
    are multihomogeneous as follows: For each monomial in the polynomial, the total degree of independent node variables coming from a given class $\alpha$ is the same as the number of nodes labeled $\alpha$ in the garden.
\end{observation}

\subsection{Composition of hanging garden diagrams}

\begin{figure}

\begin{tikzpicture}[
  font=\normalsize,
  conduit/.style={line width=1.2pt, line cap=round, line join=round},
  garden/.style={draw, fill=black!25},
  input/.style={draw, fill=white, minimum height=0.60cm, inner sep=1pt},
  nodebox/.style={draw, rounded corners=8pt, fill=black!10,
                  minimum height=0.70cm, inner sep=0pt},
  dottedframe/.style={draw, dotted, line width=0.9pt}
]
  % ------------------------------------------------------------
  % 1) Outer garden boundary/background (Type I)
  % ------------------------------------------------------------
  \draw[garden] (0.20,0.15) rectangle (9.20,11.10);

  % ------------------------------------------------------------
  % Row y-levels (top to bottom)
  % ------------------------------------------------------------
  \def\yOne{10.20}  % theta
  \def\yTwo{9.15}   % zeta
  \def\yThree{8.1} % top alpha
  \def\yFour{6.55}  % inner gamma (moved down to add space from row 3)
  \def\yFive{5.50}  % inner beta + inner input
  \def\ySix{4.50}   % inner bottom input
  \def\ySeven{3.15} % delta + outer gamma
  \def\yEight{2.05} % alpha + input + input + beta
  \def\yNine{0.95}  % final alpha

  % ------------------------------------------------------------
  % Column x-positions
  % ------------------------------------------------------------
  \def\xL{1.60}   % left branch (theta1 -> delta)
  \def\xM{4.70}   % middle (theta2 -> zeta)
  \def\xR{7.05}   % right branch (theta3 -> outer gamma) moved left

  % Zeta lower conduits
  \def\xZ2{5.40}  % zeta2 -> top alpha
  \def\xIG{3.95}  % zeta1 -> inner gamma

  % Inner region subcolumns
  \def\xIGL{3.25} % inner beta + bottom input
  \def\xIGR{4.85} % inner right input

  % Row 8 columns
  \def\xIleft{3.20}   % input under delta2
  \def\xIright{5.95}  % input under outer gamma1 (shifted left with gamma)
  \def\xB{7.85}       % beta (row 8) and final alpha (row 9), shifted left

  % ------------------------------------------------------------
  % Component centers
  % ------------------------------------------------------------
  \coordinate (TH) at (\xM,\yOne);
  \coordinate (ZE) at (\xM,\yTwo);
  \coordinate (A0) at (\xZ2,\yThree);

  \coordinate (gIn)  at (\xIG,\yFour);
  \coordinate (bIn)  at (\xIGL,\yFive);
  \coordinate (iInR) at (\xIGR,\yFive);
  \coordinate (iInB) at (\xIGL,\ySix);

  \coordinate (DE) at (\xL+0.8,\ySeven);
  \coordinate (GO) at (\xR-0.2,\ySeven);

  \coordinate (A1) at (\xL,\yEight);
  \coordinate (I1) at (\xIleft,\yEight);
  \coordinate (I2) at (\xIright,\yEight);
  \coordinate (B2) at (\xB,\yEight);

  \coordinate (A2) at (\xB,\yNine);

  % ------------------------------------------------------------
  % 2) Conduits
  % ------------------------------------------------------------

  % Outer Type I top conduit exiting the garden
  \draw[conduit] ($(TH)+(0,0.35)$) -- (\xM,11.55);

  % Theta's three lower conduits
  \draw[conduit] (\xL,\yOne-0.35) -- (\xL,\ySeven+0.00);
  \draw[conduit] (\xM,\yOne-0.35) -- (\xM,\yTwo+0.00);
  \draw[conduit] (\xR,\yOne-0.35) -- (\xR,\ySeven+0.00);

  % Zeta's two lower conduits
  \draw[conduit] (\xIG,\yTwo-0.35) -- (\xIG,\yFour+0.00);
  \draw[conduit] (\xZ2,\yTwo-0.35) -- (\xZ2,\yThree+0.00);

  % Inner-garden phenomenon: short upward stub above inner gamma
  \draw[conduit] (\xIG,\yFour+0.35) -- (\xIG,\yFour+0.85);

  % Inner gamma two lower conduits
  \draw[conduit] (\xIGL,\yFour-0.35) -- (\xIGL,\yFive+0.00);
  \draw[conduit] (\xIGR,\yFour-0.35) -- (\xIGR,\yFive-0.05);

  % Inner beta one lower conduit
  \draw[conduit] (\xIGL,\yFive-0.35) -- (\xIGL,\ySix-0.05);

  % Delta two lower conduits
  \draw[conduit] (\xL,\ySeven-0.35) -- (\xL,\yEight+0.00);
  \draw[conduit] (\xIleft,\ySeven-0.35) -- (\xIleft,\yEight-0.05);

  % Outer gamma two lower conduits (shifted left)
  \draw[conduit] (\xIright,\ySeven-0.35) -- (\xIright,\yEight-0.05);
  \draw[conduit] (\xB,\ySeven-0.35) -- (\xB,\yEight+0.00);

  % Beta (row 8) lower conduit
  \draw[conduit] (\xB,\yEight-0.35) -- (\xB,\yNine+0.00);

  % ------------------------------------------------------------
  % 3) Components
  % ------------------------------------------------------------

  \node[nodebox, minimum width=8.00cm] at (TH) {$\theta$};
  \node[nodebox, minimum width=3.00cm] at (ZE) {$\zeta$};
  \node[nodebox, minimum width=1.50cm] at (A0) {$\alpha$};

  \node[nodebox, minimum width=2.80cm] at (gIn) {$\gamma$};
  \node[nodebox, minimum width=1.50cm] at (bIn) {$\beta$};
  \node[input,  minimum width=0.90cm] at (iInR) {}; % unlabeled
  \node[input,  minimum width=0.90cm] at (iInB) {}; % unlabeled

  % Row 7: widened delta and gamma
  \node[nodebox, minimum width=3.40cm] at (DE) {$\delta$};
  \node[nodebox, minimum width=3.40cm] at (GO) {$\gamma$};

  \node[nodebox, minimum width=1.50cm] at (A1) {$\alpha$};
  \node[input,  minimum width=0.95cm] at (I1) {}; % unlabeled
  \node[input,  minimum width=0.95cm] at (I2) {}; % unlabeled
  \node[nodebox, minimum width=1.50cm] at (B2) {$\beta$};

  \node[nodebox, minimum width=1.50cm] at (A2) {$\alpha$};

  % ------------------------------------------------------------
  % 4) Dotted rectangle highlighting the inner-garden region
  % ------------------------------------------------------------
  \draw[dottedframe] (2.15,3.85) rectangle (5.8,7.4);

\end{tikzpicture}

\caption{A composition of hanging garden diagrams}
\label{fig:composition}
\end{figure}
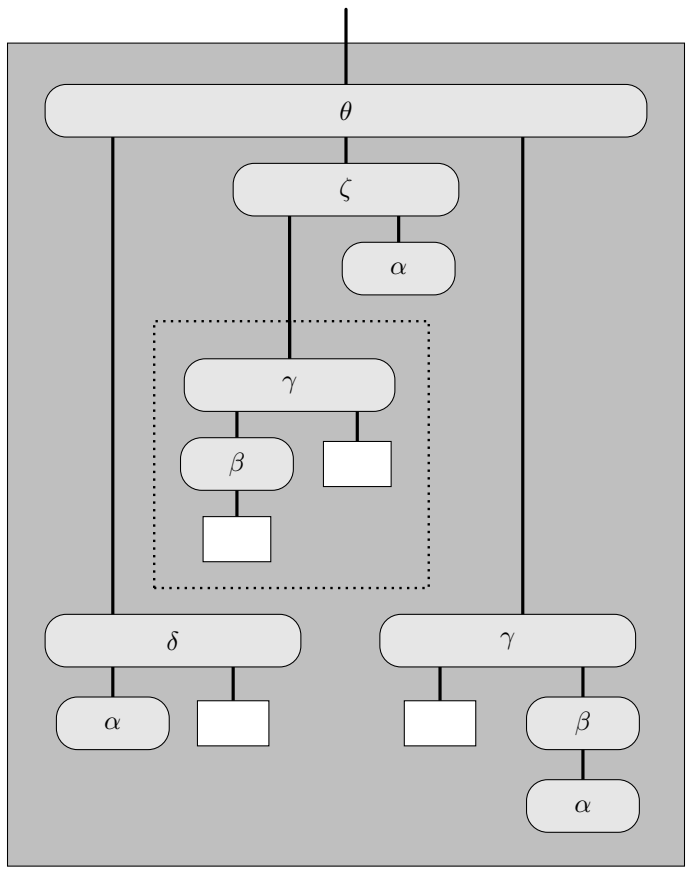

A hanging garden diagram of Type I has the right shape, including its uppermost conduit, to be placed inside one of the input boxes of another hanging garden diagram (of either Type~I or Type~II).
If the boundary of the original input box is then erased, the result is a new garden, called the \emph{composition} of the original outer and inner gardens.
\Cref{fig:composition} depicts such a composition, with a dotted line depicting the boundary of the original input box into which an inner garden was placed.

\begin{observation}
    \label{obs:compose}
    Composition of gardens is composition of functions.
\end{observation}

\begin{example}
    In \cref{fig:composition}, the outer garden of the composition, with three input boxes, represents a function $f(\mathbf{x}_1, \mathbf{x}_2, \mathbf{x}_3)$.
    Its input boxes are ordered $\mathcal{I}_1, \mathcal{I}_2, \mathcal{I}_3$, with $\mathcal{I}_2$ hosting the inner garden.
    The inner garden of the composition, with two input boxes, represents a function $g(\mathbf{x}_1, \mathbf{x}_2)$.
    The composite garden has four input boxes and represents a function
    \[
    h(\mathbf{x}_1, \mathbf{x}_2, \mathbf{x}_3, \mathbf{x}_4) = f(\mathbf{x}_1, g(\mathbf{x}_2, \mathbf{x}_3), \mathbf{x}_4).
    \]
\end{example}

Composition gives hanging gardens the structure of an operad.
This operad is reminiscent of (and was partially inspired by) the planar algebras of
Vaughan~F.~R. Jones \cite{J1999}, although
nothing like the full power of planar algebras is
required in this simple setting, and the
construction is correspondingly simpler.
Missing in gardens, for example, are loops, parity considerations with shaded and unshaded regions,
and non-trivial planar isotopies.
Uniform LSS requires only a garden-variety operad.

Any hanging garden diagram can be constructed, under composition, out of simpler gardens, down to gardens at most one of whose internal components is a node.

\subsection{Simple building blocks}

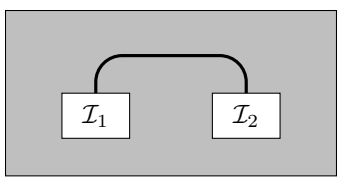
\begin{figure}
% Hanging garden diagram — Type II (Diagram 2: two-input garden)
% Minimal Type II: a single cap-shaped conduit connecting two inputs

\begin{tikzpicture}[
  font=\normalsize,
  conduit/.style={line width=1.2pt, line cap=round, line join=round},
  garden/.style={draw, fill=black!25},
  input/.style={draw, fill=white, minimum height=0.60cm, inner sep=1pt}
]
  % 1) Garden boundary/background
  \draw[garden] (0.20,0.90) rectangle (4.60,3.10);

  % Component centers
  \coordinate (I1) at (1.40,1.70); % left input
  \coordinate (I2) at (3.40,1.70); % right input

  % 2) Cap conduit (Type II), generous bend radius
  \draw[conduit, rounded corners=10pt]
    ($(I1)+(0,0.30)$) -- ($(I1)+(0,0.80)$)
    -- ($(I2)+(0,0.80)$) -- ($(I2)+(0,0.30)$);

  % 3) Inputs (draw last)
  \node[input, minimum width=0.90cm] at (I1) {$\mathcal{I}_1$};
  \node[input, minimum width=0.90cm] at (I2) {$\mathcal{I}_2$};

\end{tikzpicture}
    \caption{The unique hanging garden diagram of Type~II without internal nodes represents the inner product function.}
    \label{fig:inner-product}
\end{figure}

\subsubsection{Zero nodes, Type~II}
The simplest garden of Type II, as illustrated in \cref{fig:inner-product}, has two input boxes, no internal nodes, and a single, uppermost conduit.

\begin{observation}
    \label{obs:inner}
    The function represented by a hanging garden diagram of Type~II with no internal nodes is the inner product,
    \[
    g(\mathbf{x}_1, \mathbf{x_2}) = \mathbf{x}_1^\top \mathbf{x}_2^{\phantom{\top}}\!.
    %\left< \mathbf{x}_1, \mathbf{x}_2 \right>.
    \]
\end{observation}

\subsubsection{Zero nodes, Type~I}
The simplest garden of Type I has one input box, no nodes, and a single, uppermost conduit.
It represents the identity function,
\[
g(\mathbf{x}) = \mathbf{x},
\]
which leaves any other garden unchanged when composed with it as either inner or outer garden of the composition.
%gardens using it either as the inner or the outer garden yields as result the other garden, unchanged.

\subsubsection{One node, Type~II}
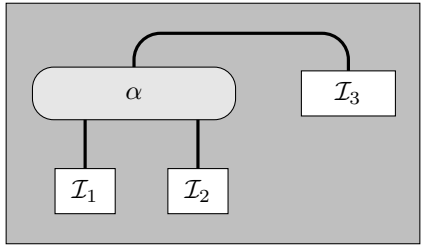
\begin{figure}
% Hanging garden diagram — Type II ("alpha")

\begin{tikzpicture}[
  font=\normalsize,
  conduit/.style={line width=1.2pt, line cap=round, line join=round},
  garden/.style={draw, fill=black!25},
  input/.style={draw, fill=white, minimum height=0.60cm, inner sep=1pt},
  nodebox/.style={draw, rounded corners=8pt, fill=black!10,
                  minimum width=2.70cm, minimum height=0.70cm, inner sep=0pt}
]
  % 1) Garden boundary/background
  \draw[garden] (0.20,0.35) rectangle (5.70,3.55);

  % Component centers (used for conduit routing)
  \coordinate (A)  at (1.90,2.35); % gamma node center
  \coordinate (B1) at (1.25,1.05); % bottom-left input center
  \coordinate (B2) at (2.75,1.05); % bottom-right input center
  \coordinate (TR) at (4.75,2.35); % top-right input center

  % 2) Conduits (slightly over-long to meet; boxes drawn later hide any intrusion)
  % Cap conduit: come down far enough to reach into the top input
  \draw[conduit, rounded corners=10pt]
    ($(A)+(0,0.35)$) -- ($(A)+(0,0.80)$)
    -- ($(TR)+(0,0.80)$) -- ($(TR)+(0,0.00)$);

  % Lower vertical conduits: extend up into the node and down into inputs
  \draw[conduit] ($(B1)+(0,-0.05)$) -- ($(B1)+(0,1.35)$);
  \draw[conduit] ($(B2)+(0,-0.05)$) -- ($(B2)+(0,1.35)$);

  % 3) Components (draw last)
  \node[nodebox] at (A) {$\alpha$};

  % Inputs numbered left-to-right across the whole diagram
  \node[input, minimum width=0.80cm] at (B1) {$\mathcal{I}_1$};
  \node[input, minimum width=0.80cm] at (B2) {$\mathcal{I}_2$};
  \node[input, minimum width=1.25cm] at (TR) {$\mathcal{I}_3$};

\end{tikzpicture}

\caption{The garden $g_\alpha$ for $\alpha$ of arity $2 + 1$.}
\label{fig:garden-alpha}

\end{figure}

Given a class $\alpha$ of node variables, with arity $k + 1$, there are two ways to form a hanging garden diagram of Type~II with only this internal node, depending on whether the node is placed to the left or to the right of the uppermost conduit.
It is convenient to place $\alpha$ on the left; call the resulting garden~$g_\alpha$.
\Cref{fig:garden-alpha} depicts $g_\alpha$ for $\alpha$ of arity $2 + 1$.
With $\alpha$ on the left, the fixed order of conduits $(1, \dots, k; k+1)$
around $\alpha$ matches the ordering of the $q = k + 1$ input boxes
$(\mathcal{I}_1, \dots, \mathcal{I}_{k+1})$.
By abuse of notation, denote the function represented by $g_\alpha$ also as
\[
g_{\alpha}(\mathbf{x}_1, \dots, \mathbf{x}_k, \mathbf{x}_{k+1}).
\]

\begin{proposition}
    \label{prop:alternating}
    For any class $\alpha$ of node variables, the function
    \[
        g_\alpha : \left( \mathbb{R}^\ad \right)^{k + 1} \longrightarrow \mathbb{R}
    \]
    is an alternating multilinear function.
\end{proposition}
\begin{proof}
    Multilinearity is by \cref{obs:multilinear}.
    Consider two matchings $A$ and $B$ of the garden $g_\alpha$ that differ only by exchanging the input vectors for two input boxes $\mathcal{I}_i$ and $\mathcal{I}_j$.
    Given an indexing of $A$, consider the indexing of $B$ that exchanges the index labels of conduits $i$ and $j$.
    Since $\alpha$ is an alternating function on its indices, and since for both $A$ and $B$ the same vector entries are extracted, the valuation of the indexed $B$ is the negative of the valuation of the indexed $A$.
    The summation over all possible ways to index the conduits of $A$ is also, applying the same index swap to $B$ in each case, a summation over all possible ways to index the conduits of $B$.
    It follows that the scalar value represented by $B$ is the negative of the scalar value represented by $A$, and that $g_\alpha$ is an alternating multilinear function as claimed.
\end{proof}

\begin{remark}
    The unique garden of Type~II with a single internal node $\alpha$ on the right (rather than left) of the uppermost conduit represents the same function $g_\alpha$ when $k$ is even, and represents the function $-g_\alpha$ (which is also an alternating multilinear function) when $k$ is odd.
\end{remark}

\begin{corollary}
    \label{cor:depend}
    If there is any linear dependency among the $k + 1$ inputs to $g_\alpha$, the value of the function is $0$.
\end{corollary}

\subsubsection{One node, Type~I}
Given a class $\alpha$ of node variables, with arity $k + 1$, there is a unique way to form a hanging garden diagram of Type~I with only this internal node.
As a function, it takes $k$ inputs from $\mathbb{R}^\ad$ and produces a vector also in $\mathbb{R}^\ad$.
Call both this garden and its associated function $\varphi_\alpha$.
Type~I gardens with a single node are the workhorse of constructing more complicated hanging garden diagrams under composition.
Referring again to \cref{fig:garden-alpha}, $g_\alpha$ can be viewed as a composition of gardens in which $\varphi_\alpha$ has been plugged into the left input of the inner product garden.

\begin{proposition}
    \label{prop:ortho}
    Given a $k$-tuple of vectors $W = (\mathbf{x}_1, \dots, \mathbf{x}_k)$,
    the output of $\varphi_\alpha(W)$ is orthogonal to the span of $W$.    
\end{proposition}
\begin{proof}
    Suppose that $\mathbf{y}$ is in the span of $W$.
    Then by \cref{cor:depend},
    \[
    g_\alpha(W; \mathbf{y}) = 0.
    \]
    On the other hand, by \cref{obs:inner} and \cref{obs:compose},
    \[
    g_\alpha(W; \mathbf{y}) = \varphi_\alpha(W)^\top \mathbf{y}.
    \]
\end{proof}
The stated purpose for hanging garden diagrams was to allow a parametrization of the LSS construction,
and \cref{prop:ortho} shows that indeed the workhorse of gardens can be used to choose a vector that is orthogonal to some set of previously chosen vectors.

\begin{example}
    In $3$-dimensional space $\ad=3$, let $\alpha$ be a node variable class of arity $2 + 1$,
    with a unique independent free variable coming from $\binom{3}{2 + 1} = 1$.
    Then $\varphi_\alpha$ is a multiple of the cross product,
    \[
    \varphi_\alpha(\mathbf{x}_1, \mathbf{x_2}) = \alpha_{1, 2; 3}\left(\mathbf{x}_1 \times \mathbf{x}_2\right).
    \]
\end{example}
Another part of the stated goal is that the parametrization should be in terms of polynomials over a free set of commuting variables, namely the independent node variables.

\begin{observation}
    If each of the inputs to $\varphi_\alpha$ is a vector whose entries are integer polynomials over independent node variables, then the output is also a vector whose entries are integer polynomials over independent node variables, including those of~$\alpha$.
    The overall multidegree of the output vector polynomials is the sum of the multidegrees of the input vector polynomials, plus degree $1$ in the independent node variables of $\alpha$.
\end{observation}

Finally, the goal is not just to parametrize LSS but in particular to parametrize \emph{uniform} LSS, which requires $\varphi(v_j) = \mathbf{0}$ whenever the previously chosen vectors are dependent.

\begin{proposition}
    Suppose that the vectors in a $k$-tuple $W = (\mathbf{x}_1, \dots, \mathbf{x}_k)$
    are linearly dependent.
    Then $\varphi_\alpha(W) = \mathbf{0} \in \mathbb{R}^\ad$.
\end{proposition}
\begin{proof}
    Let $\mathbf{y} \in\mathbb{R}^\ad$ be any vector at all.
    By \cref{cor:depend}, \cref{obs:inner}, and \cref{obs:compose} as above,
    \[
    \varphi_\alpha(W)^\top \mathbf{y} = 0.
    \]
    But the only vector orthogonal to every $\mathbf{y}$ is the zero vector
    $\varphi_\alpha(W) = \mathbf{0} \in \mathbb{R}^\ad$.
\end{proof}

The necessary tools are now in place to discuss hanging garden diagrams as a parametrization of uniform LSS.

\section{LSS gardens, greedy orderings, and non-cancellation}
\label{sec:LSS-garden}
\subsection{Only LSS gardens}
A bit of refreshed context first, to set the stage:

Recall from \cref{subsec:parametrization} the stated goal to parametrize uniform LSS by obtaining entries of the representation matrix $R$ and the Gram matrix $T = [t_{ij}]  = R^\top \!R$ as polynomials over free commuting variables.

Recall from \cref{subsec:LSS} the notation
\[
W_j = (w_1, w_2, \dots, w_k) := (v_i \in V \ |\  i < j \mbox{ and } v_i \not\sim v_j ),
\]
with $k(v_j) = |W_j|.$

Recall from \cref{subsec:nodevar} that
each vertex $v_1, v_2, v_3, \dots$ of $G$ should be associated to a specific class $\alpha, \beta, \gamma, \dots$ of node variables.
To make the association more precise, denote the node variable class associated to $v_j$ as $\alpha(v_j)$, yielding for example
\[
\alpha(v_1) = \alpha;\  \ \alpha(v_2) = \beta, \ \ \alpha(v_3) = \gamma,
\]
and so forth.
The arity of $\alpha(v_j)$ is $k(v_j) + 1$.

\begin{example}
    Given $k(v_6) = 2$ with $\ad=4$, the free independent node variables associated to $v_6$, of arity $2 + 1$, can be written
    \[
\alpha(v_6)_{1, 2; 3}, \ \alpha(v_6)_{1, 2; 4}, \ \alpha(v_6)_{1, 3; 4}, \mbox{ and } \alpha(v_6)_{2, 3; 4}.
\]
\end{example}

From this point forward, the only hanging gardens considered will be those that specifically model uniform LSS.
\begin{definition}
    An \emph{LSS garden} is a hanging garden diagram subject to the following constraints:
\begin{itemize}
    \item There are no input boxes.
    \item For Type~I, the top node is labeled by the node class of some vertex of $G$.
    \item For Type~II, the two top nodes are both labeled by the node classes of vertices of $G$.
    \item For any node $A$ in the garden with label $\alpha(v_j)$, consider the corresponding ordered list
    $W_j = (w_1, \dots, w_k)$ of non-adjacent precursors to $v_j$.
    Then the lower conduits of $A$, of which there are $k = k(v_j)$ in number, must be connected {\bf in order} to nodes labeled
    \[
    \alpha(w_1), \dots, \alpha(w_k).
    \]
\end{itemize}
\end{definition}
These constraints give exact recursive specifications for only
\begin{itemize}
    \item $n$ gardens of Type~I,
    \begin{itemize}
        \item representing the vectors $\varphi(v_j)$ chosen during uniform LSS,
        \item each of which is column $j$ of the representation matrix $R$,
    \end{itemize}
    \item and $n^2$ gardens of Type~II,
    \begin{itemize}
        \item which are pairwise inner products of LSS gardens of Type~I,
        \item and are the entries $t_{ij}$ of the Gram matrix $T = R^\top \! R$.
    \end{itemize}
\end{itemize}
Since $v_i$ can belong to $W_j$ only when $v_i$ strictly precedes $v_j$, the recursion does terminate,
but with repeated branching it can result in a number of nodes in the overall diagram that is exponential in $n$.

\begin{proposition}
    \label{prop:zero-case}
    Let $g$ be an LSS garden of Type~II with top node labels $\alpha(v_j)$ and $\alpha(v_i)$, $i < j$, such that $v_j$ and $v_i$ are non-adjacent in $G$.
    Then $g$ evaluates to the scalar $0$.
\end{proposition}
\begin{proof}
    This follows by the alternating nature of nodes, \cref{prop:alternating}:
    Call the top node labeled $\alpha(v_i)$ and its entire downstream $A$.
    There is also a node labeled $\alpha(v_i)$ somewhere directly below the top node labeled $\alpha(v_j)$; call this lower copy and its entire downstream $B$.
    Then $A$ and $B$ are identical, represent the same vector, and when composed into an alternating function return $0$.

    This can also be seen more directly: for every indexing of the conduits of $g$ for which the conduit labels above $A$ and above $B$ are different, consider the permutation that exchanges those two conduit labels and also exchanges all the labels in the lower parts of $A$ with all the corresponding labels in the lower parts of $B$.
    This permutation gives a bijection between all the ways to index conduits that result in a nonzero monomial, and for each such pair the valuation monomials are identical except for a change in sign from the exchanged conduit labels adjacent to the node labeled $\alpha(v_j)$.
    In summation, the monomials cancel in pairs and give an overall scalar value of $0$,
    as expected for entry $t_{ij}$ of the Gram matrix of an orthogonal representation of $G$ with vectors $v_j$ and $v_i$ non-adjacent.
\end{proof}

\subsection{What can go wrong?}
\begin{figure}
% Hanging garden diagram for monomial cancellation example

\begin{tikzpicture}[
  font=\normalsize,
  conduit/.style={line width=1.2pt, line cap=round, line join=round},
  garden/.style={draw, fill=black!25},
  nodebox/.style={draw, rounded corners=8pt, fill=black!10,
                  minimum width=2.10cm, minimum height=0.70cm, inner sep=0pt}
]
  % 1) Garden boundary/background (wider, with more top margin)
  \draw[garden] (0.20,0.35) rectangle (6.40,7.25);

  % Column x-positions (more space between columns)
  \coordinate (L) at (2.10,0); % left column
  \coordinate (R) at (4.60,0); % right column

  % Node centers (left column)
  \coordinate (dL) at ($(L)+(0,6.00)$);
  \coordinate (gL) at ($(L)+(0,4.50)$);
  \coordinate (bL) at ($(L)+(0,3.00)$);
  \coordinate (aL) at ($(L)+(0,1.50)$);

  % Node center (right column alpha aligned with delta)
  \coordinate (aR) at ($(R)+(0,6.00)$);

  % 2) Conduits (draw before nodes)
  % Top Type II cap conduit, with generous headroom
  \draw[conduit, rounded corners=10pt]
    ($(dL)+(0,0.35)$) -- ($(dL)+(0,0.9)$)
    -- ($(aR)+(0,0.9)$) -- ($(aR)+(0,0.00)$);

  % Left column vertical conduits
  \draw[conduit] ($(gL)+(0,0.00)$) -- ($(dL)+(0,-0.35)$);
  \draw[conduit] ($(bL)+(0,0.00)$) -- ($(gL)+(0,-0.35)$);
  \draw[conduit] ($(aL)+(0,0.00)$) -- ($(bL)+(0,-0.35)$);

  % 3) Nodes (draw last)
  \node[nodebox] at (dL) {$\delta$};
  \node[nodebox] at (gL) {$\gamma$};
  \node[nodebox] at (bL) {$\beta$};
  \node[nodebox] at (aL) {$\alpha$};

  \node[nodebox] at (aR) {$\alpha$};

\end{tikzpicture}

\caption{A hanging garden diagram for \cref{ex:P4} and \cref{ex:P4-garden}.}
\label{fig:garden-p4}
\end{figure}
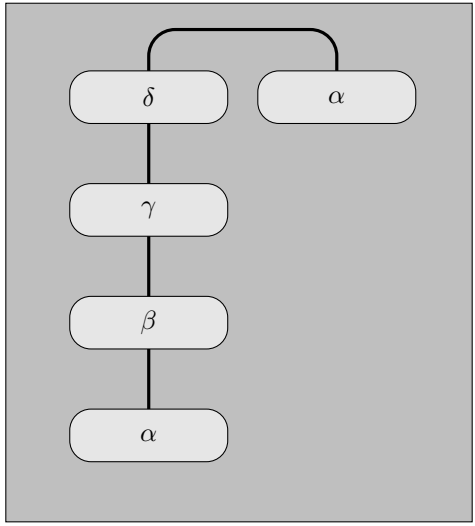
Having accomplished a full parametrization of uniform LSS, it is instructive to revisit \cref{ex:P4}, in which uniform LSS does give an orthogonal representation, and even succeeds weakly, but does not succeed strongly.
This will illustrate, from the garden and polynomial perspective, the sort of problem that must somehow be avoided if a faithful orthogonal representation in codimension $\gdel(G)$ is to be guaranteed.

\begin{example}
\label{ex:P4-garden}
As before, the graph $G$ considered is the complement of $P_4$, which is another copy of $P_4$ but with an unusual vertex ordering that allows
each vector choice $\varphi(v_j)$ to be constrained only by orthogonality with the immediately previous vector.
The failure happens in dimension $\ad = 2$, and
the offending entry of the Gram matrix is $t_{41}$,
which should represent an edge of $G$, but which LSS forces to have value~$0$.
This phenomenon should be detected by the hanging garden corresponding to entry $t_{41}$, as depicted in \cref{fig:garden-p4}.
Observe that no node has more than a single lower conduit, which is why uniform LSS succeeds at least weakly in dimension $\ad = 2$.
In particular, it is possible to index the conduits of this zero-input garden in some ways that yield a nonzero monomial.
There are, to be precise, exactly two ways in which this can be done.
If the uppermost conduit is labeled by index $1$, then all other conduit indices must alternate, yielding the valuation
\[
\alpha_1\,\delta_{2; 1}\, \gamma_{1; 2}\, \beta_{2; 1}\, \alpha_2
= \alpha_1\, \alpha_2\, \beta_{1;2}\, \gamma_{1, 2}\, \delta_{1; 2}.
\]
If on the other hand the uppermost conduit is labeled by index $2$, then again all other conduit indices must alternate, yielding the valuation
\[
\alpha_2\,\delta_{1; 2}\, \gamma_{2; 1}\, \beta_{1; 2}\, \alpha_1
= - \alpha_1\, \alpha_2\, \beta_{1;2}\, \gamma_{1, 2}\, \delta_{1; 2}.
\]
The only two nonzero monomials cancel each other, 
yielding $t_{41} = 0$ overall.
\end{example}
It is tempting to suppose that the above example is small enough that the cancellation is ``accidental'' and would not be expected in sufficiently large or complicated examples.
There do exist much more complicated examples, however---arising, for example, from essentially any incidence theorem for finite arrangements in linear geometry---where a much larger set of much larger non-zero monomials still conspires to exactly cancel in pairs.
In order to guarantee a faithful orthogonal representation, there must be some mechanism for avoiding altogether this troublesome pairwise annihilation.

As already alluded to, the saving mechanism will come from a combination of
\begin{itemize}
    \item greedy vertex orderings, and
    \item uniqueness of the leading term in an appropriate monomial order.
\end{itemize}
\subsection{Greedy orderings and anti-greedy orderings}
\label{subsec:antigreedy}
From this point forward, it will be assumed that the vertex ordering $V(G) = (v_1, \dots, v_n)$ is greedy.

Instead of talking about a greedy ordering on $G$ per se, it is sometimes useful to reason instead about an ``anti-greedy'' ordering on its complement $H = \overline{G}$, particularly since the conduit connections in the gardens being considered typically come from non-adjacencies in $G$, corresponding to the edges of $H$.
Conceptually, an anti-greedy ordering of $H$ introduces as few new edges from $H$ as possible each time a new vertex $v_j$ is introduced.
To make this precise, it is useful to reiterate \cref{def:local}---but for $H$ instead of $G$, and with new notation---and then to adapt \cref{def:greedy} to its opposite, using the new notation.
\begin{definition}
\label{def:anti-local}
Given an ordering
$V(H) =(v_1, \dots, v_{n})$
and integers $i$ and $j$ such that $1 \le i \le j \le n$,
the \emph{local degree in $H$ of $v_j$ at stage $i$}, denoted $k_i(v_j)$, is the number of
edges in $H$ between $H\left[\left\{v_1, \dots, v_{i - 1}\right\}\right]$ and $v_j$,
or in other words the size of the set $\{v_hv_j \in E(H) \st h < i\}$.
The local degree in $H$ of $v_j$ at its own stage $j$, or $k_j(v_j)$, is abbreviated to $k(v_j)$.
\end{definition}

\begin{remark}
    The notation $k(v_j)$ appears to clash with an existing notation related to the graph $G=\overline{H}$, but in fact this is the identical function
    \[k(v_j) = |W_j|.\]
\end{remark}

\begin{definition}
An ordering $V(H) =(v_1, \dots, v_{n})$ is called \emph{anti-greedy}
if, for all $1 \le i \le j \le n$, the weak inequality
\[
k(v_i) \le k_i(v_j)
\]
holds.
\end{definition}

\begin{observation}
    Given a simple graph $G$ with complement $H = \overline{G}$, an ordering of the common vertex set $V(G) = V(H)$ is a greedy ordering of $G$ if and only if it is an anti-greedy ordering of $H$.
\end{observation}

\begin{observation}
    \label{obs:kmono}
    For $i < j$,
    \[
    k_i(v_j) \le k_j(v_j) = k(v_j)
    \]
    always holds.
    If in addition the vertex ordering is greedy for~$G$ (and thus anti-greedy for~$H$),
    \[
    k(v_i) \le k_i(v_j) \le k(v_j)
    \]
    also holds, and the function $k(v_j)$ is weakly increasing in $j$.
\end{observation}

The monotonicity of $k(v_j)$ for greedy orderings is one reason why only the (maximum possible) degree of the final vertex is what defines $\gdel(G)$, in determining a dimension sufficiently large for LSS to succeed strongly.

\subsection{Greedy structure in LSS gardens}

The defining anti-greedy inequality
gives a crucial constraint on the structure of a greedy LSS garden.

\begin{proposition}
    \label{prop:w_s}
    Suppose that the vertices of $G$ are in greedy order, and for any vertex $v_j$ with $k = k(v_j)$ consider the ordered list
    \[
    W_j = (w_1, \dots, w_k)
    \]
    of non-adjacent precursors to $v_j$.
    Then for any $s \in \{1, \dots, k\}$,
    the inequality
    \[
    k(w_s) < s
    \]
    holds strictly.
\end{proposition}
\begin{proof}
    Let $i$ be such that $w_s = v_i$, implying that
    \[k(w_s) = k(v_i) \le k_i(v_j)\]
    must hold, as the defining inequality making the ordering on $H = \overline{G}$ anti-greedy.
    The list $W_j$ is carefully ordered in precisely such a way that
    \[
    k_i(v_j) = s - 1,
    \]
    giving
    \[
    k(w_s) < s
    \]
    as claimed.
\end{proof}

\begin{corollary}
    \label{cor:barren}
    Within any LSS garden for a greedily-ordered graph $G$,
    consider any node labeled $\alpha(v_j)$ and its
    $k = k(v_j)$ ordered children.
    For any $s \in \{1, \dots, k\}$,
    child number $s$ has strictly fewer than $s$ children of its own.
\end{corollary}

\begin{proposition}
    \label{prop:w_s_plus}
    Suppose that the vertices of $G$ are in greedy order, and for any vertex $v_j$ with $k = k(v_j)$ consider the ordered list
    \[
    W_j = (w_1, \dots, w_k)
    \]
    of non-adjacent precursors to $v_j$.
    Let $v_i$ be a vertex adjacent to $v_j$ in $G$ (therefore not occurring in $W_j$), and let $s$ be such that
    exactly $s - 1$ of the vertices in the list $W_j$ precede $v_i$ in the vertex ordering of $G$.
    Then
    \[
    k(v_i) < s
    \]
    holds strictly.
\end{proposition}
\begin{proof}
    As in the previous proof,
    \[k(v_i) \le k_i(v_j)\]
    must hold, and as before
    the list $W_j$ is ordered such that
    \[
    k_i(v_j) = s - 1,
    \]
    giving
    \[
    k(v_i) < s
    \]
    as claimed.
\end{proof}

\begin{corollary}
    \label{cor:barren_plus}
    Let $g$ be an LSS garden of Type~II for a greedily-ordered graph $G$ whose
    top nodes are labeled $\alpha(v_j)$ and $\alpha(v_i)$, $i < j$,
    with $v_j \sim v_i$ in $G$.
    Consider the top node labeled $\alpha(v_j)$ and its
    $k = k(v_j)$ ordered children,
    with $W_j = (w_1, \dots, w_k)$.
    Let $s$ be the position that $v_i$ would occupy if it were inserted into the list $W_j$ in order.
    Then the top node labeled $\alpha(v_i)$ has fewer than $s$ children.
\end{corollary}

\subsection{A lexicographic term order}
In order to show that not all monomials in a garden polynomial cancel in pairs,
it is useful to single out one particular monomial that has no canceling partner,
which is done by specifying a term ordering according to which the special monomial has strictly highest priority.
The general principal is that early vertices and then small numbers have priority, lexicographically.
That is to say, the ordering is lexicographic, both by
node variable class
\[
\alpha(v_1) > \alpha(v_2) > \dots > \alpha(v_n),
\]
and, within each node variable class, by lexicographic ordering of the independent node variables, with lower indices bestowing higher priority.
For example:
\[
\alpha(v_3)_9 > \alpha(v_4)_{1,7;8} > \alpha(v_4)_{2, 3; 6} > \alpha(v_4)_{2, 4; 5} > \alpha(v_5)_{1, 2; 3}.
\]

As it happens (by \cref{obs:multihom}), all monomials in a given garden polynomial have the same total degree within each node variable class, which degree therefore does not affect priority.
The import of the lexicographic ordering is that the ordering within later node classes, when comparing monomials, matters only if the product of independent node variables within earlier classes is identical on either side.

\begin{lemma}[Unique Monomial Lemma]
    Let $g$ be a greedy LSS garden of Type~II with top node labels $\alpha(v_j)$ and $\alpha(v_i)$, $i \le j$, such that either $v_i = v_j$ or $v_j \sim v_i$ in $G$,
    and let $\ad$ be at least $k(v_n) + 1$.
    Then among all index labelings of the conduits of $g$, the lexicographically highest priority nonzero monomial occurs uniquely, such that no other nonzero monomial can annihilate it.
\end{lemma}
\begin{proof}
    By \cref{obs:kmono}, if $\ad \ge k(v_n) + 1$ holds in a greedy ordering, then $\ad \ge k(v_j) + 1$ holds for all $j \le n$, and there are sufficient indices available that, at least a priori, no node is forced to have two adjacent conduits with the same index label.
    (Equivalently, though not necessary for the proof of the lemma, with $|W_j| < \ad$ for all $v_j$, the subspace $S_j$ from which $\varphi(v_j)$ is chosen will always have dimension at least $1$, and uniform LSS will succeed at least weakly.)
    
    Some small cases will be examined before considering a full inductive proof---both to seek a pattern, and to establish an inductive base case.
    Due to the lexicographic ordering first by node variable class, no indexing that happens in the conduits surrounding nodes labeled $\alpha(v_2) = \beta$ matters until all conduits surrounding nodes labeled $\alpha(v_1) = \alpha$ have been indexed, and so forth.
    The first vertex $v_1$ has $k(v_1) = 0$ and all nodes labeled $\alpha$ have arity $0 + 1$, or a single conduit above.
    For the highest priority, all such conduits must be labeled $1$, if this can be done in a way that yields a nonzero monomial---and it can, because $\alpha$ never labels two nodes that are direct siblings of each other.
    
    If $v_1$ happens to be an isolated vertex in $G$, then $k(v_2) = 1$, every node labeled $\beta$ has a conduit labeled $1$ below it, leading down to a node labeled $\alpha$, every list of children begins with a first node labeled $\alpha$ that has a conduit labeled $1$ above it, and everywhere that a node labeled $\beta$ occurs, it will be the second node in the list of children and the conduit above it will, out of lexicographic priority, but avoiding two conduits labeled $1$ adjacent to $\beta$ and avoiding two conduits labeled $1$ adjacent to the node above $\beta$, be labeled $2$. Since there are never two sibling nodes both labeled $\beta$, no index greater than $2$ is yet necessary.
    
    In most cases, though, where $v_1$ is not an isolated vertex in $G$, in a greedy ordering $v_2$ must be adjacent to $v_1$, so that $v_2$ also has $k(v_2) = 0$ with arity $0 + 1$.
    By lexicographic priority, all nodes labeled $\beta$ must be surrounded fully by already-indexed conduits before any consideration is given to nodes labeled $\gamma$.
    The nodes labeled $\beta$ would most prefer to have their unique upper adjoining conduit labeled by index $1$, if this can be done while allowing a nonzero monomial.
    Suppose, though, that some node, say labeled $\alpha(v_4) = \delta$, has below it both a node labeled $\alpha$ and a node labeled $\beta$.
    They cannot both have their upper conduits given index $1$, or this $\delta$ will have two adjoining conduits with the same index label, yielding a zero monomial for the indexed valuation.
    The first child labeled $\alpha$ took priority and already got $1$ for its conduit index, so the second child labeled $\beta$ must take its second-choice index label $2$ on its upper conduit.
    On the other hand, anywhere that a node labeled $\beta$ is the first child, it can be granted its first preference and be given upper conduit index label $1$.

    The pattern so far is that below any node the lower conduits to $(w_1, \dots, w_k)$ always get labels $1$ through $k$ in order.
    This is what will be proven inductively, with one exception:
    Consider the higher- (or equal-) valued of the two top nodes, labeled $\alpha(v_j)$, with $k = k(v_j)$ and $W_j = (w_1, \dots, w_k)$.
    The label $\alpha(v_i)$ of the other top node, $i \le j$, does not occur immediately below among the lower labels $\alpha(w_s)$ since $v_i$ is adjacent (or equal) to $v_j$ in $G$.
    Create a list
    \[W_j^\prime = (w_1^\prime, \dots, w_{k+1}^\prime)\]
    in which $v_i$ has been inserted into $W_j$, in its proper order, with $v_i = w_t^\prime$ for some $t$ in the range $1 \le t \le k + 1$.
    Then what is to be proven, in addition to the non-exceptional cases,
    is that those lower conduits of the first top node (labeled $\alpha(v_j)$) from $1$ through $t - 1$ get index labels $1$ through $t - 1$ in order,
    that the uppermost conduit from $\alpha(v_j)$ to $\alpha(v_i)$ gets index label $t$, and that lower conduits $t$ through $k$ get index labels $t + 1$ through $k + 1$ in order.
    In other words, the conduits adjacent to the top $\alpha(v_j)$ get indexed $1$ through $k + 1$ in order of the priority of the node labels at the other side of the conduit.
    (As a sanity check, if there were a tie in priority, with $v_i \not \sim v_j$ in $G$ and therefore $v_i$ equal to some $w_s$, the two ways of settling the tie would lead to two opposite-signed and annihilating top-priority monomials.)
    In the case where $v_i = v_j$, it is ambiguous which of the two sets of lower conduits the above discussion was supposed to apply to, but in that case there is in fact perfect symmetry between them: On both sides, the children get conduits labeled $1$ to $k$ in order, and the uppermost $\cap$-shaped conduit between them gets label $k + 1$.

    Assume by way of induction that up to but not including some vertex number $v_u$, $1 < u \le n$, conduit indices have been assigned strictly according to lexicographic priority, and that the result so far has always matched the desired pattern.
    Now for lexicographic priority, all conduits adjacent to nodes labeled $\alpha(v_u)$ must be assigned indices, before nodes labeled $\alpha(v_{u+1})$ or higher start to matter.
    By induction, all nodes with labels lower than $\alpha(v_u)$ have had both their lower conduits and their upper conduit indexed, conforming to the pattern.
    By induction, every node with label $\alpha(v_u)$ has all its lower conduits but not its upper conduit indexed.
    By induction, wherever a node labeled $\alpha(v_u)$ occurs as a child (or possibly as a top node or pseudo-child, if $u = i$),
    all preceding children have had their upper conduits indexed in child order.
    If $v_u = v_j$, either the top node labeled $\alpha(v_j)$ already has all adjacent conduits indexed, in which case the induction is complete, or $v_u = v_j = v_i$ and after labeling the uppermost conduit by $k + 1$, the induction is complete.
    Consider any one node labeled $\alpha(v_u)$; call it $C$ and call its parent $P$---or, in the case that $u=i$ and $C$ is a top node, give the name $P$ to its pseudo-parent, the other top node labeled $\alpha(v_j)$.
    The lowest available index for the conduit above $C$ that would not give identical index labels to two conduits adjacent to $P$ is the index corresponding to the numerical position of $C$ among the children of $P$.
    By either \cref{cor:barren} or \cref{cor:barren_plus}, this choice of conduit index above $C$ does not give a repeated index adjacent to $C$, either.
    These are exactly the lexicographically highest priority choices that can safely be made for all nodes labeled $\alpha(v_u)$, and all such choices preserve the desired pattern.
    Lexicographically forced choices thus continue the desired pattern through all stages of induction.

    The choices made were unique in every case; no other indexing of the conduits can achieve a monomial of such high priority while also avoiding two conduits adjacent to the same node with the same index.
\end{proof}

\begin{remark}
    \label{rem:signs}
    The list $W_j^\prime$ in the exceptional case of the proven structure is exactly the list $W_{j; i}$ from \cref{subsec:strong-success}, but put in order, and since all other node variables of the proven structure are already in order, the sign of the permutation that puts $W_{j; i}$ in order determines the sign of the dominating monomial.
\end{remark}

\subsection{Proof of main results}
\label{subsec:mainproof}

All necessary pieces are now in place to prove the promised main result first stated in \cref{subsec:mainstate}.
To recapitulate:

\mainthm*

\begin{proof}
    The Unique Monomial Lemma requires $d \ge k(v_n) + 1$, which is satisfied by the hypotheses of the Main Theorem with
    \[
    k(v_n) + 1 = (n - 1) - \deg (v_n) + 1 = (n - 1) - (n - \ad) + 1 = \ad.
    \]
    (This does mean, for example, that $S_n$ will have dimension only $1$ and the last vector choice $\varphi(v_n)$ will be unique up to scaling.)

    The purpose of hanging garden diagrams, which parametrize uniform LSS by polynomials, is to show that for a greedy ordering, uniform LSS in a dimension $d$ of at least $n - \deg(v_n)$ produces a faithful orthogonal representation with probability $1$ and succeeds strongly.
    The fact that it does give an orthogonal representation---or in other words, that for $v_j \not \sim v_i$ the scalar value of a hanging garden $g$ of Type~II is $0$---is taken care of by \cref{prop:zero-case}.
    The remaining cases involve a hanging garden diagram $g$ of Type~II whose top nodes are labeled $\alpha(v_j)$ and $\alpha(v_i)$ with
    \begin{itemize}
        \item $v_i = v_j$, where showing that each value $g$ is non-zero with probability $1$ will show that each $\varphi(v_j)$ is a non-zero vector with probability $1$, establishing that uniform LSS succeeds weakly, and
        \item $v_i \sim v_j$ in $G$, where showing that each value $g$ is non-zero with probability $1$ will show that the Gram matrix has nonzero entries everywhere that it should, establishing that uniform LSS produces a faithful orthogonal representation with probability $1$.
    \end{itemize}
    Establishing these two things implies, by \cref{prop:strongly},
    that uniform LSS succeeds strongly, as required for the main results.
    In both cases, the Unique Monomial Lemma establishes that the scalar value of $g$, which is an integer polynomial over independent node variables, is not the zero polynomial.
    Uniform LSS chooses $\varphi(v_j)$ from each orthogonal subspace $S_j$ according to some reasonable probability distribution, making generic choices that with probability $1$ result in a non-zero value for every non-zero polynomial.
    For each above choice of hanging garden $g$, and simultaneously for all of them, with probability $1$ the scalar values of the hanging gardens are non-zero, giving nonzero entries $t_{ij}$ exactly where needed to demonstrate that uniform LSS under these conditions succeeds strongly, as required.
\end{proof}

\section{Generalizations and concluding remarks}
\label{sec:more}
\subsection{Numerical LSS}
Although hanging gardens can be useful for determining when LSS should succeed,
the polynomials of exponential size that they model are unwieldy for producing actual examples.
When LSS does work, it just works, and concrete examples can easily be generated numerically, either using floating-point (if approximate examples suffice) or rational arithmetic.

\subsection{Indefinite orthogonal representations}
The version of uniform LSS considered for the main results here always produces a positive semidefinite Gram matrix,
which in combinatorial matrix theory is generally (but not quite universally) the most difficult sort of matrix to produce with a given nullity for a given pattern~$G$.
Given an indefinite form $E$ that is diagonal with diagonal entries $\epsilon_i = \pm1$,
the analysis of hanging gardens can be modified to track monomial signs coming from the~$\epsilon_i$,
allowing one to model the polynomials of a Gram matrix with a chosen inertia other
than positive semidefinite.

\subsection{Other infinite fields}
In \cite{AIM2006} it is supposed that the Delta Conjecture ought to extend to other infinite fields.
The hanging garden methodology should apply equally well to such extensions, since it produces non-zero polynomials, which in any infinite field evaluate generically to a non-zero value.

\subsection{Hermitian gardens}
Hanging garden diagrams can be used to model an LSS process that produces a Hermitian rather than a real symmetric matrix, with the caveat that each conduit within a diagram should be oriented to reflect the difference between an entry and its complex conjugate. All conduits around a given node should be oriented either toward it or away from it.
Since the connectivity within a hanging garden is just a tree (rather than the more complicated topologies allowed, for example, in planar algebras), this is always possible.

\subsection{Particular sign patterns}
Given any lexicographic monomial ordering and given any finite set of polynomials, it is possible to produce parameter values of widely different orders of magnitude in such a way that the highest-priority monomial of each polynomial dominates all other terms, meaning in particular that the sign of that one monomial determines the overall sign of the evaluated polynomial.
As mentioned in \cref{rem:signs},
the sign of the dominating term from the Unique Monomial Lemma is the sign of the permutation $P$
that puts the list $W_{j; i}$ in order.
That sign can be determined by examining the pattern of the Gram matrix $T$:
Find the corresponding non-zero entry $t_{ij}$ above the diagonal;
then the number of swaps to perform the permutation $P$ is the same as the number of zeros in column $j$ below $t_{ij}$ and above the diagonal.
The Unique Monomial Lemma can thus give not only a faithful Gram matrix, but a faithful Gram matrix with this particular sign pattern---or any sign pattern that can be formed under some choice of the signs of those independent node variables that appear in dominating monomials.

\subsection{Dominating monomials for non-greedy orderings}
The proof of the Unique Monomial Lemma depends crucially on a greedy ordering of $G$,
but a similar process can be carried out for LSS gardens modeling a non-greedy ordering---it will just fail in low dimensions.
For the garden in \cref{fig:garden-p4}, for example, the process would go as follows:
\begin{itemize}
    \item Both conduits above nodes labeled $\alpha$ get index $1$.
    \item The conduit above $\beta$ gets index $2$.
    \item But now the conduit above $\gamma$ and below $\delta$ can get neither index $1$ nor index $2$, requiring $\ad \ge 3$.
\end{itemize}
The process does yield a unique dominating monomial, but only in a dimension that is strictly higher than the naive bound
\[
d \ge \max \{k(v_i) + 1\}
\]
that suffices when the vertex ordering is greedy.

For a given graph $G$ with a given ordering $V=(v_1, \dots, v_n)$ that is not necessarily greedy,
what is the lowest dimension that always allows unique dominating monomials according to this process---or some other process?
Is there a better way (lower dimension bound) to ensure that strong LSS succeeds for non-greedy orderings?
For every graph, ordering, and dimension, if uniform LSS does succeed strongly, that fact may be much easier to detect with numerical simulation than by trying to find a unique dominating monomial, and in fact uniform LSS might sometimes strongly succeed in dimensions low enough that the highest-priority monomial is not unique and does get annihilated, as long as some lower-priority monomial does not.

\subsection{Degrees of genericity}
The results of LSS show that a Gram matrix $T$ with pattern $G$ whose columns are fully generic (that is, such that any set of columns up to cardinality the rank of $T$ is independent) can only be produced with nullity up to $\kappa(G)$, rather than $\delta(G)$ or $\gdel(G)$.
Weaker than full genericity is the genericity of SAP, which has been used to great effect in combinatorial matrix theory.
Add to these another generic condition, upper-column generic, which is sufficiently weak to allow nullity $\gdel(G)$ and yet is still stronger than SAP,
\[
\gdel(G) \le \mathrm{M}_+^u(G) \le \nu.
\]
Little is known about just how much stronger than SAP it can be, or what sort of graphs give gaps of what sizes in the associated parameter values. 
Also, while $\nu(G)$ is automatically minor-monotone,
the same is likely not true of $\mathrm{M}_+^u(G)$.

\appendix

\section{Approximate greedegree is NP-hard}
\label{sec:NP}

The claim is that calculation of $\gdel(G)$, even to within any fixed approximation, is NP-hard.
The question is first formulated as a decision problem by specifying a threshold degree $\delta$ and considering the challenge of producing a greedy ordering in which all vertices of low degree $\le \delta$ are chosen before the final vertex.
The exact value of $\gdel(G)$ can then be calculated using logarithmically many repetitions of the decision problem.

\newcommand{\gd}{\normalfont\textsc{greedegree}}
\newcommand{\fac}{\normalfont\textsc{f}}
\newcommand{\gdf}{\normalfont\textsc{greedegree-factor(\fac)}}
\newcommand{\gde}{\normalfont\textsc{greedegree-fraction(\ensuremath{\varepsilon}})}
\newcommand{\pair}{\normalfont\textsc{pair-avoid}}

\vskip 0.15in
\noindent
\textbf{Problem:} $\gd$

\noindent \textbf{Input:} A simple graph $G=(V, E)$ and a threshold degree $\delta$.

\noindent \textbf{Output:}
\begin{itemize}
    \item YES if some greedy ordering of the vertices of $G$ ends with degree $> d$.
    \item NO if all greedy orderings of $G$ end with degree $\le \delta$.
\end{itemize}

The question of calculating $\gdel(G)$ to within some wide approximation is formulated as a promise problem in two ways:

\vskip 0.15in
\noindent
\textbf{Problem:} $\gdf$ (for fixed real approximation factor $\fac \ge 1$)

\noindent \textbf{Input:} A simple graph $G=(V, E)$ and a threshold degree $\delta$.

\noindent \textbf{Promise:} An input satisfies the promise if one of these holds:
\[
    \gdel(G) > \fac \delta \ \ \ \mbox{ or } \ \ \  \gdel(G) \le \delta.
\]

\noindent \textbf{Output:}
\begin{itemize}
    \item YES, if some greedy ordering of the vertices of $G$ ends with a vertex of degree strictly greater than $\fac \delta$.
    \item NO, if every greedy ordering of $G$ ends with a vertex of degree at most $\delta$.
    \item unspecified, if the promise does not hold.
\end{itemize}

\vskip 0.15in
\noindent
\textbf{Problem:} $\gde$ (for fixed real $0 < \varepsilon \le \frac12$)

\noindent \textbf{Input:} A simple graph $G=(V, E)$.

\noindent \textbf{Promise:} An input satisfies the promise if one of these holds:
\[
    \frac{\gdel(G)}{|V|} > 1 - \varepsilon \ \ \ \mbox{ or } \ \ \  \frac{\gdel(G)}{|V|} \le \varepsilon.
\]

\noindent \textbf{Output:}
\begin{itemize}
    \item YES, if some greedy ordering of the vertices of $G$ ends with a vertex of degree strictly greater than $(1 - \varepsilon)|V|$.
    \item NO, if every greedy ordering of $G$ ends with a vertex of degree at most $\varepsilon|V|$.
    \item unspecified, if the promise does not hold.
\end{itemize}

As was mentioned in \cref{subsec:greedegree}, for $n \ge 3$ any pair of adjacent vertices can both be avoided as the final vertex of a greedy ordering simply by choosing them as the first and second vertices.
This is not at all the case, however, for a pair of non-adjacent vertices, which gives the hard case for the next decision problem.

\vskip 0.15in
\noindent
\textbf{Problem:} $\pair$

\noindent \textbf{Input:} A simple graph $G=(V, E)$ and a pair of vertices $a, b \in V$.

\noindent \textbf{Output:}
\begin{itemize}
    \item YES, if some greedy ordering of the vertices of $G$ ends with a vertex other than $a$ or $b$.
    \item NO, if every greedy ordering of $G$ either ends with $a$ or ends with $b$.
\end{itemize}

\begin{theorem}
\label{thm:nphard}
The decision problem $\gd$ is NP-complete.
The decision problem $\pair$ is NP-complete.
For any real $\fac  \ge 1$, the promise problem $\gdf$ is NP-hard.
For any real $0 < \varepsilon \le \frac{1}{2}$, the promise problem $\gde$ is NP-hard.
\end{theorem}
The proof of \cref{thm:nphard} that concludes this Appendix will make use of a formula graph whose construction, together with some of its key properties, will first be explained.

\newcommand{\cla}{k}
\newcommand{\ato}{t}
\newcommand{\fil}{f}

\begin{definition}
    \label{def:construct}
A \emph{formula graph} is a simple graph constructed as follows:

\noindent \textbf{Input:} A non-trivial Boolean formula $\Phi$ in conjunctive normal form, together with either
an approximation fraction $0 < \varepsilon \le \frac12$ or
an approximation factor $\fac  \ge 1$.

\noindent \textbf{Output:} A formula graph $G(\Phi, \varepsilon)$ (possibly with $\varepsilon$ in terms of $\fac $), a lower threshold degree $\delta$, and an upper threshold degree $\Delta$.

The number of vertices in the constructed graph will be denoted $n$.
For a requested approximation fraction $\varepsilon$,
the lower threshold should satisfy 
\[ \delta \le \varepsilon n \]
and the upper threshold should satisfy
\[ \Delta \ge (1 - \varepsilon) n.\]
In the case that an approximation factor $\fac  \ge 1$ is instead requested, this is first converted to an approximation fraction
\[
\varepsilon = \frac{1}{\fac  + 1},
\]
which satisfies $0 < \varepsilon \le \frac12$
and which by the above will ensure $\Delta \ge \fac \delta$.

Let $\cla$ be the count of clauses in $\Phi$, and number the clauses
\[
    \Phi = c_1 \wedge c_2 \wedge \dots \wedge c_\cla.
\]
Each clause $c_i$ is a disjunction of literals, say $w_i$ of them (for ``width''), and each literal $\ell$ is either a positive atom $x_i$ or the negation $\neg x_i$ of an atom.
Let $\ato$ be the count of atoms that occur anywhere, so that the complete set of possible literals is
\[ x_1, \neg{x_1},\,  x_2, \neg{x_2}, \, \dots,\,  x_\ato, \neg{x_\ato}.\]
Without loss of generality, assume that no literal occurs more than once within the same clause.
For each literal $\ell$, let $p(\ell)$ (for ``prevalence'') represent the number of clauses in which $\ell$ occurs.
Let $m$ represent the maximum value attained by any width $w_i$ or doubled prevalence $2p(\ell)$:
\[
   m = \max \{w_1, \dots, w_\cla, 2p(x_1), 2p(\neg x_1), \dots, 2p(x_\ato), 2p(\neg x_\ato)\}.
\]

Some number $\fil$ of filler vertices will be required, with the total number of vertices given by
\[
n = n(\fil) = 2 + 2\ato + 2\cla + \fil.
\]
The filler count $\fil$ will depend on the formula $\Phi$ and the approximation fraction $\varepsilon$; it is chosen as the smallest integer such that all of these bounds are satisfied:
\begin{eqnarray*}
    \fil &\ge& m; \\
    n(\fil) &\ge& \frac1\varepsilon(2\ato + \cla); \\
    n(\fil) &\ge& \frac1\varepsilon(m + 3).
\end{eqnarray*}
The first inequality, which is independent of $\varepsilon$, will play a role in the proof of \cref{thm:nphard}.
The second and third inequalities will relate to the degree thresholds $\delta$ and $\Delta$, respectively.

The $n$ vertices of $G(\Phi, \varepsilon)$ are named as follows:
\begin{itemize}
\item (2 vertices) A left critical vertex $a$ and a right critical vertex $b$.
\item ($2\ato$ vertices) For each literal $x_i$ or $\neg x_i$, a vertex that by abuse of notation is also named $x_i$ or $\neg{x_i}$, respectively.
\item ($2\cla$ vertices) For each clause $c_i$, a left clause vertex $\alpha_i$ and a right clause vertex $\beta_i$.
\item ($\fil$ vertices) Filler vertices $y_1, \dots, y_{\fil}$.
\end{itemize}
The graph $G(\Phi, \varepsilon)$ is quite dense and so it is convenient to specify it by the adjacencies of its complement $H(\Phi, \varepsilon)$, as follows:
\begin{itemize}
\item The two critical vertices $a$ and $b$ are adjacent in $H$.
\item For each atom $x_i$, the positive literal $x_i$ and the negative literal $\neg x_i$ are adjacent in $H$.
\item For each clause $c_i$ of the form
\[
    c_i = \ell_1 \vee \dots \vee \ell_{w_i},
\]
\begin{itemize}
    \item the left clause vertex $\alpha_i$ is adjacent in $H$ to the left critical vertex $a$ and to all of the literals $\ell_1, \dots, \ell_{w_i}$; similarly,
\end{itemize}
\item Each filler vertex $y_i$ is adjacent in $H$ to both critical vertices $a$ and $b$.
\end{itemize}
Note that $m$ is was chosen in precisely such a way that $m + 1$ is the maximum degree obtained by any left clause vertex, right clause vertex, or literal in $H$.

It remains to output a lower threshold degree
\[
    \delta = 2\ato + \cla
\]
and an upper threshold degree
\[
    \Delta = n - m - 3.
\]
\end{definition}
\begin{proposition}
    \label{prop:Gdegrees}
    The degree thresholds of $G(\Phi, \varepsilon)$ satisfy
    \[
        \delta \le \varepsilon n \ \ \mbox{ and } \ \ \Delta \ge (1 - \varepsilon) n
    \]
    and the degree thresholds of $G(\Phi, \frac{1}{\fac + 1})$ 
    satisfy 
    \[
    \Delta \ge \fac \delta.
    \]
    Furthermore, every vertex in $G(\Phi, \varepsilon)$ either has degree no more than $\delta$ or has degree strictly larger than $\Delta$.
\end{proposition}
\begin{proof}
The bounds on $\delta$ and $\Delta$ are a direct consequence of requirements
\[
    n(\fil) \ge \frac1\varepsilon(2\ato + \cla) = \frac1\varepsilon \delta
\]
and
\[
    n(\fil) \ge \frac1\varepsilon(m + 3) = \frac1\varepsilon(n(\fil) - \Delta)
\]
on the number $\fil$ of filler vertices.

It remains to show that vertex degrees are indeed constrained by the lower and upper thresholds.
Tabulating non-adjacencies in $H(\Phi, \varepsilon)$ yields a count for the
degree of each vertex in $G(\Phi, \varepsilon)$ as follows:
\begin{itemize}
    \item Each critical vertex ($a$ or $b$) has degree \[n - \fil - \cla - 2 \ = \ 2\ato + \cla = \delta.\]
    \item Each literal $\ell$ has degree \[n - 2p_\ell - 2 \ > \  n - m - 3 = \Delta.\]
    \item Each clause vertex $\alpha_i$ or $\beta_i$ has degree \[n - w_i - 2\  > \ n - m - 3 = \Delta.\]
    \item Each filler vertex $f_i$ has degree \[n - 3 \ > \ n - m - 3= \Delta.\]
\end{itemize}

\end{proof}

Since $\gdel(G)$ is always the degree of some vertex of $G$,
\cref{prop:Gdegrees} has consequences for the promise problems.
\begin{corollary}
    \label{cor:promise}
    For $\gdf$,
    the pair $(G(\Phi, \frac{1}{\fac + 1}), \delta)$ fulfills the promise.
\end{corollary}
\begin{corollary}
    \label{cor:fraction}
    For $\gde$,
    the graph $G(\Phi, \varepsilon)$ fulfills the promise.
\end{corollary}

An example will serve to illustrate.

\begin{example}
    \label{ex:construct}
    Take for $\Phi$ the Boolean formula in conjunctive normal form
    \[
    \Phi = (x_1 \vee x_2 \vee x_3) \wedge (\neg x_1 \vee x_2) \wedge
        (x_1 \vee \neg x_3) \wedge (\neg x_2 \vee x_3) \wedge
        (\neg x_1 \vee \neg x_2 \vee \neg x_3),
    \]
    a small non-trivial example of a formula with no satisfying assignment.
    (The cyclic implications of the middle three clauses enforce the logical equivalence of all three atoms, but by the first and last clauses there must be both a true atom and a false atom among them.)
    This gives $\cla=5$, $w_i \in \{2, 3\}$, $\ato=3$, $p(\ell) \in \{2\}$, and $m=\max \{2, 3, 2(2)\} =4$.
    Take the trivial approximation factor $\fac = 1$, giving
    \[\varepsilon = \frac{1}{\fac +1} = \frac12.\]
    The filler vertex count $\fil$ should be the lowest integer satisfying
    \begin{eqnarray*}
        \fil & \ge& m = 4; \\
        n(\fil) = 18 + \fil &\ge& \frac{1}{\varepsilon} (2\ato + \cla) =
            2 (11); \\
        n(\fil) = 18 + \fil & \ge& \frac{1}{\varepsilon} (m + 3) =
            2 (7).
    \end{eqnarray*}
    The first two constraints coincide and give
    $\fil = 4$ filler vertices with $n=22$ vertices in total.
    The vertices of $G(\Phi, 1/2)$ are as follows:
    \begin{itemize}
        \item $a$, $b$,
        \item $x_1, \neg x_1, \, x_2, \neg x_2, \, x_3, \neg x_3$,
        \item $\alpha_1, \beta_1, \alpha_2, \beta_2, \alpha_3, \beta_3,
        \alpha_4, \beta_4, \alpha_5, \beta_5$,
        \item $y_1, y_2, y_3,$ and $ y_4.$
    \end{itemize}
    The low degree threshold and high degree threshold are $\delta=2\ato + \cla = 11$ and $\Delta = n - m - 3 = 15$, respectively, reflecting the fact that $a$ and $b$ have degree $11$ and the smallest other degree is $16$.
    
    The complement graph $H(\Phi, 1/2)$ on $22$ vertices is depicted by \cref{fig:h-Phi-F}.

\begin{figure}
        \caption{The complement graph $H(\Phi, 1/2)$ for \cref{ex:construct}.}
        \label{fig:h-Phi-F}
\begin{tikzpicture}[
    scale=1.0,
    vertex/.style={circle, draw, inner sep=1.5pt},
    >=stealth
]

%--- a and b (midline, medium width) ---
\node[vertex,label=above left:{$a$}] (a) at (-5,0) {};
\node[vertex,label=above right:{$b$}] (b) at ( 5,0) {};

%--- literals (vertical column below) ---
% within pair: 0.5 ; between pairs: 1.0
\node[vertex,label={[yshift=0.5pt]above:{$x_1$}}]      (x1)  at (0,-0.8) {};
\node[vertex,label={[yshift=-0.8pt]below:{$\neg x_1$}}] (nx1) at (0,-1.4) {};
\node[vertex,label={[yshift=0.5pt]above:{$x_2$}}]      (x2)  at (0,-2.6) {};
\node[vertex,label={[yshift=-0.5pt]below:{$\neg x_2$}}] (nx2) at (0,-3.2) {};
\node[vertex,label={[yshift=0.8pt]above:{$x_3$}}]      (x3)  at (0,-4.4) {};
\node[vertex,label={[yshift=-0.5pt]below:{$\neg x_3$}}] (nx3) at (0,-5.0) {};

% literal pair edges
\draw (x1) -- (nx1);
\draw (x2) -- (nx2);
\draw (x3) -- (nx3);

%--- slanted alpha column (left, below) ---
% from just below x1 to just above ¬x3
\node[vertex,label={[xshift=-3pt,yshift=-8pt]center:{$\alpha_1$}}] (a1) at (-2.8,-0.8) {};
\node[vertex,label={[xshift=-3pt,yshift=-8pt]center:{$\alpha_2$}}] (a2) at (-3.15,-1.7) {};
\node[vertex,label={[xshift=-3pt,yshift=-8pt]center:{$\alpha_3$}}] (a3) at (-3.5,-2.6) {};
\node[vertex,label={[xshift=-3pt,yshift=-8pt]center:{$\alpha_4$}}] (a4) at (-3.85,-3.5) {};
\node[vertex,label={[xshift=-3pt,yshift=-8pt]center:{$\alpha_5$}}] (a5) at (-4.2,-4.4) {};

%--- slanted beta column (right, mirror of alpha) ---
\node[vertex,label={[xshift=3pt,yshift=-8pt]center:{$\beta_1$}}] (b1) at ( 2.8,-0.8) {};
\node[vertex,label={[xshift=3pt,yshift=-8pt]center:{$\beta_2$}}] (b2) at ( 3.15,-1.7) {};
\node[vertex,label={[xshift=3pt,yshift=-8pt]center:{$\beta_3$}}] (b3) at ( 3.5,-2.6) {};
\node[vertex,label={[xshift=3pt,yshift=-8pt]center:{$\beta_4$}}] (b4) at ( 3.85,-3.5) {};
\node[vertex,label={[xshift=3pt,yshift=-8pt]center:{$\beta_5$}}] (b5) at ( 4.2,-4.4) {};

%--- f-vertices (vertical column above midline) ---
% spacing d = 0.7; f6 at y = 2d = 1.4
%\node[vertex,label={[yshift=-2pt]above:{$f_1$}}] (f1) at (0,4.8) {};
%\node[vertex,label={[yshift=-2pt]above:{$f_1$}}] (f1) at (0,4.0) {};
\node[vertex,label={[yshift=-2pt]above:{$y_1$}}] (y1) at (0,3.2) {};
\node[vertex,label={[yshift=-2pt]above:{$y_2$}}] (y2) at (0,2.4) {};
\node[vertex,label={[yshift=-2pt]above:{$y_3$}}] (y3) at (0,1.6) {};
\node[vertex,label={[yshift=-2pt]above:{$y_4$}}] (y4) at (0,0.8) {};

%--- base edge between a and b ---
\draw (a) -- (b);

%--- edges from f_i to a and b ---
\foreach \i in {1,...,4}{
  \draw (a) -- (y\i);
  \draw (y\i) -- (b);
}

%--- edges from a to each alpha_i ---
\foreach \i in {1,...,5}{
  \draw (a) -- (a\i);
}

%--- edges from b to each beta_i ---
\foreach \i in {1,...,5}{
  \draw (b) -- (b\i);
}

%--- adjacency: alpha_1, beta_1 to x_1, x_2, x_3 ---
\foreach \v in {a1,b1}{
  \draw (\v) -- (x1);
  \draw (\v) -- (x2);
  \draw (\v) -- (x3);
}

%--- adjacency: alpha_2, beta_2 to ¬x_1 and x_2 ---
\foreach \v in {a2,b2}{
  \draw (\v) -- (nx1);
  \draw (\v) -- (x2);
}

%--- adjacency: alpha_3, beta_3 to x_1 and ¬x_3 ---
\foreach \v in {a3,b3}{
  \draw (\v) -- (x1);
  \draw (\v) -- (nx3);
}

%--- adjacency: alpha_4, beta_4 to ¬x_2 and x_3 ---
\foreach \v in {a4,b4}{
  \draw (\v) -- (nx2);
  \draw (\v) -- (x3);
}

%--- adjacency: alpha_5, beta_5 to ¬x_1, ¬x_2, ¬x_3 ---
\foreach \v in {a5,b5}{
  \draw (\v) -- (nx1);
  \draw (\v) -- (nx2);
  \draw (\v) -- (nx3);
}

\end{tikzpicture}

    \end{figure}
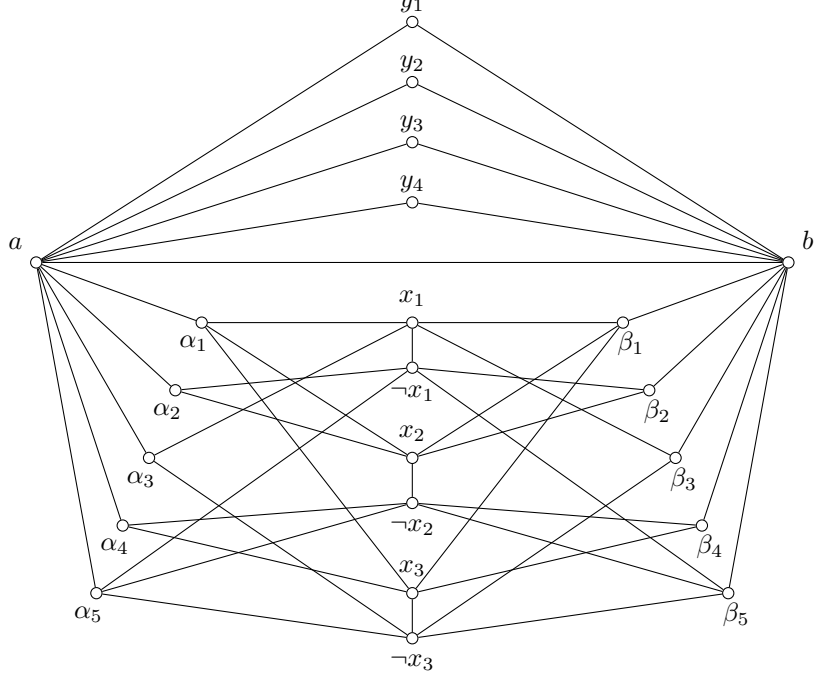
\end{example}

Just as the sparse complement $H = \overline{G}$ is easier to depict than the dense graph $G$,
it will be easier to reason about anti-greedy orderings on $H$ than to reason about greedy orderings on $G$.
Anti-greedy orderings are described in \cref{subsec:antigreedy}.

\begin{observation}
    \label{obs:anti-greedy}
    Any greedy ordering of a graph $G$ begins with some locally maximal clique.
    Any anti-greedy ordering of a graph $H$ begins with some locally maximal independent set.
\end{observation}

\begin{proof}[Proof of \cref{thm:nphard}]
To show that $\gd$ and $\pair$ belong to NP requires a certificate for each.
This is easily provided by the vertex ordering itself,
%Taking as certificate the vertex ordering,
which can be checked in polynomial time to be a greedy ordering and, for $\gd$, to have final degree $> \delta$ or, for $\pair$, to end with a vertex other than $a$ or $b$.
The proof that the promise problems $\gdf$ and $\gde$ are NP-hard, and hence that $\gd$ is NP-complete, is by reduction from \textsc{sat},
using a construction that also achieves reduction of $\pair$ from \textsc{sat},
showing that in addition $\pair$ is NP-complete.
Given any non-trivial Boolean formula $\Phi$ in conjunctive normal form, the construction of
\cref{def:construct} provides a graph $G = G(\Phi, \varepsilon)$ on $n$ vertices, with complement $H = H(\Phi, \varepsilon)$, and threshold degrees $\delta$ and $\Delta$
such that $G$ has exactly two vertices $a$ and $b$ of low degree $\delta \le \varepsilon n$, all other vertices having degree strictly greater than $\Delta \ge \fac\delta = (1 - \varepsilon) n$.
By \cref{cor:promise},
the pair $(G, \delta)$ does satisfy the promise required for $\gdf$, and by \cref{cor:fraction} the graph $G$ similarly does satisfy the promise required for $\gde$.
The claim is that for each promise problem, and for $\pair$ with $G$ and the pair $(a, b)$, the answer is YES if and only if $\Phi$ has a satisfying assignment.

The challenge presented in each case, for the answer YES, is to find a greedy vertex ordering of $G$ such that both of the low-degree vertices $a$ and $b$ are chosen before the last vertex in the ordering.
Equivalently, in the complement $H$, the challenge is to find an anti-greedy vertex ordering in which neither $a$ nor $b$ is chosen last, with the difficulty arising from the fact that in $H$ they are the two vertices of high degree,
and must compete with other vertices of overall low degree in order to be chosen early, in spite of an anti-greedy local policy that prioritizes low adjacency counts first.
While considering this challenge in what follows, adjacencies, degrees, and local degrees will always be considered in the context of the complement graph $H$.
Recall that in $H(\Phi, \varepsilon)$ each left clause vertex, right clause vertex, or literal has degree at most $m + 1$.

Suppose first that $\Phi$ does have a satisfying assignment $(\ell_1, \dots, \ell_\ato)$, where each literal $\ell_i$ is either $x_i$ or $\neg x_i$.
Observe that the vertices $S = \{a, \ell_1, \dots, \ell_\ato\}$ form an independent set in $H$, from which it follows that the sequence $(a, \ell_1, \dots, \ell_\ato)$ is a valid start to an anti-greedy ordering.
Can $b$ be chosen next?
It is adjacent only to $a$ in $S$,
and thus has local degree $1$.
It can be chosen next if and only if no other vertex outside of $S$ has local degree $0$.
Each right clause vertex $\beta_i$ is adjacent to at least one of the literals in the satisfying assignment.
Every other vertex in $V(H) \setminus S$ is either adjacent to $a$ or is a literal that is adjacent to its opposite literal in $S$.
It follows that any next chosen vertex would have local degree at least $1$, which means that $b$ can indeed be chosen next, rather than being postponed to the end of the ordering.
When $\Phi$ does have a satisfying assignment, the answer to each promise problem and to $\pair$ is YES.

Suppose now that $\Phi$ does not have any satisfying assignment, and consider the possibilities for an anti-greedy ordering of $H$,
with the goal of showing that one of $a$ or $b$ will be forced to be chosen last.
By \cref{obs:anti-greedy}, such an ordering must begin with a locally maximal independent set $S$, which cannot contain both of the adjacent vertices $a$ and $b$.
Without loss of generality, assume $b \not \in S$.

\vskip 5pt

\noindent \textbf{Case 1:} $a \in S$.
Since no assignment $(\ell_1, \dots, \ell_\ato)$ is a satisfying assignment,
the locally maximal independent set $S$ must contain at least one right clause vertex $\beta_i$.
The local degree of $b$, which is adjacent both to $a$ and to $\beta_i$, is at least $2$, but the local degree of every filler vertex $y_i$ is only $1$, and will remain so until $b$ is chosen.
It follows that all of the filler vertices must be chosen before $b$ is chosen,
after which the local degree of $b$ will be at least $\fil + 2$.
The first constraint on $\fil$ in the construction of $H(\Phi, \varepsilon)$ by \cref{def:construct} is
\[
    \fil \ge m,
\]
and so once all the filler vertices are by force chosen, the local degree of $b$ is at least $m + 2$.
But all other remaining vertices have total degree at most $m + 1$,
and therefore local degree at most $m + 1$ at any stage.
In the case $a \in S$, the second high-degree vertex $b$ must unavoidably be chosen last in an anti-greedy ordering,
resulting in the answer NO to either promise problem and NO to $\pair$.

\vskip 5pt

\noindent \textbf{Case 2:} $a \not \in S$ and $b \not \in S$.
Since $S$ is locally maximal, it must contain every filler vertex $y_i$.
Since every independent set of literals represents a subset of an assignment, and since no assignment is satisfying, there must be some clause $c_i$ of $\Phi$ such that both $\alpha_i$ and $\beta_i$ are in $S$.
Starting with $S$, proceed to make anti-greedy choices of vertices until one of $a$ or $b$, without loss of generality $a$, is chosen.
From that point on, the local degree of $b$---which is adjacent to every filler vertex, to $\beta_i$, and to $a$---will as before be at least $m + 2$, whereas no other remaining vertex has degree greater than $m + 1$.
As before, the second high-degree vertex $b$ must unavoidably be chosen last in an anti-greedy ordering, resulting in the answer NO to either promise problem and NO to $\pair$.

The claim is verified: For each promise problem, the promise is fulfilled, and for each promise problem and \textsc{pair-avoid}, the answer is YES if and only if $\Phi$ has a satisfying assignment.
It follows that for any $\fac \ge 1$ the promise problem $\gdf$ is NP-hard, that for any $0 < \varepsilon \le \frac12$ the promise problem $\gde$ is NP-hard, that the decision problem $\gd$ is NP-complete, and that the decision problem $\pair$ is NP-complete.
\end{proof}

% Does this translate into some tight connection between \textsc{sat} and strong success of LSS? Not in any obvious way, no.}

\bibliography{ams-greedegree}{}
\bibliographystyle{amsplain}

\end{document}